\documentclass[11pt]{amsart}
\usepackage[T1]{fontenc}
\usepackage{amssymb,mathrsfs}
\usepackage{booktabs,enumitem}
\usepackage[hidelinks]{hyperref}
\usepackage{microtype}
\usepackage{tabularx}
\usepackage{array}
\numberwithin{equation}{section}
\newtheorem{theorem}{Theorem}[section]
\newtheorem{proposition}[theorem]{Proposition}
\newtheorem{lemma}[theorem]{Lemma}
\newtheorem{corollary}[theorem]{Corollary}
\newcommand{\R}{\mathbb R}
\newcommand{\N}{\mathbb N}
\newcommand{\E}{\mathcal E}
\newcommand{\mci}{\mathcal I}
\newcommand{\LL}{\mathcal L}
\newcommand{\NN}{\mathcal N}
\newcommand{\Dspace}{\mathcal D}
\newcommand{\Wspace}{\mathcal W}
\newcommand{\dd}{\,dx}
\newcommand{\sph}{\mathrm{sph}}
\newcommand{\loc}{\mathrm{loc}}
\newcommand{\glo}{\mathrm{glo}}
\newcommand{\tr}{\operatorname{tr}}
\newcommand{\Id}{\operatorname{Id}}
\newcommand{\Ric}{\operatorname{Ric}}
\newcommand{\Rm}{\operatorname{Rm}}
\newcommand{\diag}{\operatorname{diag}}
\newcommand{\qhat}{\widehat q}
\newcommand{\ip}[2]{\langle #1,#2\rangle}
\DeclareMathOperator{\sgn}{sgn}
\setlist[itemize]{leftmargin=*,itemsep=3pt,topsep=4pt}

\makeatletter
\renewcommand*{\l@subsection}{\@tocline{2}{0pt}{3pc}{}{}}
\makeatother
\allowdisplaybreaks[4]
\title[Sharp vanishing order of the Weyl tensor]{On the sharp vanishing order of the Weyl tensor at a blow-up point for the Yamabe equation in all dimensions}
\author{Liuwei Gong}
\address{Department of Mathematics, The Chinese University of Hong Kong, Shatin, NT, Hong Kong}
\email{lwgong@math.cuhk.edu.hk}
\author{YanYan Li}
\address{Department of Mathematics, Rutgers University, 110 Frelinghuysen Road, Piscataway, NJ, USA}
\email{yyli@math.rutgers.edu}
\author{Anthony Matos}
\address{Department of Mathematics, Rutgers University, 110 Frelinghuysen Road, Piscataway, NJ, USA}
\email{amm939@math.rutgers.edu}
\date{September 15, 2026}
\hypersetup{pdftitle={Sharp vanishing order of the Weyl tensor},pdfauthor={Liuwei Gong, YanYan Li, Anthony Matos},pdfsubject={Weyl-jet blow-up on arbitrary manifolds of positive Yamabe type}}

\begin{document}
\begin{abstract}
We  determine the sharp vanishing order of the Weyl tensor at a blow-up point of solutions to the Yamabe equation in every dimension $n\ge 25$. On every closed,  connected, smooth $n$-dimensional manifold $M$ of positive Yamabe type and at any prescribed point $p\in M$, we construct a smooth metric $ \widehat g$ and a sequence of smooth positive solutions to the Yamabe equation blowing up at $p$, such that the first nonzero covariant derivative of the Weyl tensor at $p$ has any prescribed order $l$ satisfying $l\ge 6$ when $25\le n\le 29$, $l \ge 4$ when $30\le n\le 51$, and $l\ge 2$ when $n\ge 52$. Moreover, $|\nabla_{ \widehat g}^{l}W_{ \widehat g}|>0$ everywhere on $M$.

Earlier constructions realize every strictly higher order as the first nonzero order, while the known Weyl vanishing results exclude every lower first nonzero order. Thus the endpoint cases $l=6,4,2$ were the only ones left open, and are settled here.
In particular, the order-two examples show that the compactness criterion $|W|+|\nabla W|>0$ is sharp.

At the endpoints, the construction arranges global and localized metric perturbations to interact at leading order in a limiting functional. For $l=2,4$,  anisotropic polynomial ansatzes involving two orthogonal algebraic Weyl tensors remove a distinguished translation degeneracy and produce  $C^1$-stable  saddle points. For $l=6$,  a different polynomial ansatz produces a  nondegenerate local minimum. Higher prescribed orders are obtained by separating the perturbation that determines the limiting functional from the one that prescribes the Weyl derivative.
\end{abstract}
\maketitle
\setcounter{tocdepth}{2}
\enlargethispage{2pt}
\begingroup
\small
\tableofcontents
\endgroup
\section{Introduction}

Let $(M,g)$ be a closed smooth Riemannian manifold of dimension $n\ge3$. The Yamabe problem is to find a metric of constant scalar curvature in the conformal class $ [g]:=\{g_u=u^{4/(n-2)} g:0<u\in C^\infty(M)\}.$

Write $R_g$, $W_g$, and $\nabla_g$ for the scalar curvature, Weyl tensor, and covariant derivative associated with $g$. The scalar curvature of $g_u$ satisfies
\begin{equation}\label{eq:conformal-law}\
c(n):=\frac{n-2}{4(n-1)},\ \  L_g:=\Delta_g-c(n)R_g,\ \ R_{g_u}=-c(n)^{-1}u^{-\frac{n+2}{n-2}}L_gu,
\end{equation}
where $\Delta_g:=\operatorname{div}_g\nabla_g$; thus $-L_g$ is the conformal Laplacian. Positive solutions of
\begin{equation}\label{eq:yamabe}
-L_gu=n(n-2)u^{\frac{n+2}{n-2}}
\end{equation}
produce metrics with scalar curvature $4n(n-1)$. A point $p$ is a blow-up point of a sequence $v_\nu$ if $v_\nu(p_\nu)\to\infty$ for some $p_\nu\to p$. The Yamabe constant of the conformal class is
\begin{equation*}
Y(M,[g]):=\inf_{u\in H^1(M)\setminus\{0\}}\frac{c(n)^{-1}\displaystyle\int_M (|\nabla_gu|^2+c(n)R_gu^2)\,dV_g}{\left(\displaystyle\int_M|u|^{\frac{2n}{n-2}}\,dV_g\right)^{\frac{n-2}{n}}};
\end{equation*}
$Y(S^n)$ denotes the value for the round conformal class.

The existence theory was initiated by Yamabe and Trudinger and completed by Aubin and Schoen \cite{Yamabe,Trudinger,Aubin,Schoen}. The subsequent compactness problem concerns the normalized positive solutions of \eqref{eq:yamabe}. On the round sphere,  this solution space  is noncompact because the conformal diffeomorphism group is noncompact. Schoen established compactness in the locally conformally flat case when $(M,g)$ is not conformally equivalent to the round sphere, and developed the general blow-up framework \cite{SchoenLCF}. Li-Zhu, Druet, Li-Zhang, and Marques subsequently proved compactness for $n\le 11$ and, in all dimensions, under the condition $|W_g|+|\nabla W_g|>0$ \cite{LiZhu,Druet1, LiZhang1,Druet2,LiZhang2, MarquesApriori, LiZhang3}. Berti--Malchiodi constructed  noncompact examples  with  $C^k(S^n)$ metrics for $n\ge 4k+3$ \cite{BertiMalchiodi} . For smooth metrics, Brendle constructed noncompact examples for $n\ge52$ \cite{Brendle}, and Brendle--Marques extended the construction to  $n\ge25$ \cite{BrendleMarques}. Khuri-Marques-Schoen established compactness for  $n\le 24$ \cite{KMS}.

A crucial ingredient in the proofs of the compactness results in dimensions $n\le 24$, without Weyl assumptions, is the positive mass theorem. It is applied to the metric \(G_p^{4/(n-2)}g\) on \(M\setminus\{p\}\), where $G_p$ is the Green's function of the conformal Laplacian with pole $p$. This metric is scalar-flat and, by the compactness analysis, asymptotically flat with the required decay. Schoen and Yau \cite{SchoenYauMass} first proved positive mass in dimension three using minimal surfaces; their minimal hypersurface method  extends to dimensions $3\leq n\leq7$. Witten \cite{Witten} subsequently proved the spin case in all dimensions. Lesourd–Unger–Yau \cite{LesourdUngerYau} extended the positive mass theorem to manifolds which have one asymptotically flat end and possibly additional  arbitrary ends. Schoen and Yau \cite{SchoenYauSingularities} and Lohkamp \cite{Lohkamp}  proposed proofs of the positive mass theorem in arbitrary dimension without the spin assumption. The generic regularity results for area-minimizing hypersurfaces of Smale \cite{Smale}, Chodosh-Mantoulidis-Schulze \cite{ChodoshMantoulidisSchulze}, and Chodosh-Mantoulidis-Schulze-Wang \cite{ChodoshMantoulidisSchulzeWang} yield positive mass in ambient dimensions $8$ through $11$. Bi-Hao-He-Shi-Zhu \cite{BiHaoHeShiZhu} proved it through dimension $19$. More recently, Brendle-Wang \cite{BrendleWang} proved the positive mass theorem in all dimensions. The Liouville-type theorem of Caffarelli–Gidas–Spruck  \cite{CGS} is another crucial ingredient in all the  compactness results above. Blow-up constructions include Gong–Li's constant-scalar-curvature conformal metrics with unbounded volumes  \cite{GongLi}, as well as   Almaraz's scalar-flat Yamabe examples on manifolds with boundary in dimensions $n\ge25$ \cite{Almaraz}.

\subsection{Main result and sharpness}
The sharper local question asks for the lowest possible order of the first nonzero Weyl derivative at a blow-up point. Li--Zhang \cite{LiZhang2} proved that  $W_g(p)=0$ at a blow-up point $p$ when $n\ge6$, and that $|W_g(p)|=|\nabla_gW_g(p)|=0$ when $n\ge8$; Marques \cite{MarquesApriori} independently  proved the first assertion. Thus $|W|+|\nabla W|>0$ on $M$ ensures  compactness for $n\ge 8$. On the other hand, Marques constructed examples with first nonzero Weyl derivative of any order $l\ge 7$ for $25\le n\le29$, $l\ge 5$ for $30\le n\le51$, and $l\ge 3$ for $n\ge52$ \cite{MarquesBlowup}. These examples suggest that compactness can fail in every dimension $n\ge 52$ even when $|W|+|\nabla W|+|\nabla^2W|+|\nabla^3W|>0$ everywhere on $M$.

\medskip

\noindent{\bf Question.}\
{\it Does there exist an $n_0 \in \mathbb N$ such that, in every dimension $n\ge n_0$, the condition $|W|+|\nabla W|+|\nabla^2W|>0$ on $M$ ensures compactness of the normalized  solution space of \eqref{eq:yamabe}  ? }

\medskip

Our order-two examples, presented in Theorem \ref{thm:global-blowup} below,  answer this question negatively and show that the compactness criterion $|W|+|\nabla W|>0$ is sharp. As  Marques observed in \cite{MarquesBlowup}, the argument of Khuri, Marques, and Schoen \cite{KMS} forces the Weyl tensor and its covariant derivatives to vanish through order $5$ when $25\le n\le29$ and through order $3$ when $30\le n\le51$. Together with the first-order vanishing result of Li--Zhang \cite{LiZhang2}  and Marques's examples, this  left precisely the endpoints in the last column   unresolved:

\begin{center}
\small \renewcommand{\arraystretch}{1.28}
\begin{tabularx}{\textwidth}{
>{\centering\arraybackslash}p{0.15\textwidth} >{\raggedright\arraybackslash}X >{\raggedright\arraybackslash}X >{\centering\arraybackslash}p{0.16\textwidth}} \toprule \textbf{Dimension} & \textbf{Orders forced to vanish} & \textbf{Previously realized first nonzero orders} & \textbf{Endpoint realized here}\\
\midrule $25\le n\le 29$ & $0, 1,2,3,4,5$ & Every $l\ge 7$ & $l=6$ \\
\addlinespace[0.25em] $30\le n\le 51$ & $0,1,2,3$ & Every $l\ge 5$ & $l=4$ \\
\addlinespace[0.25em] $n\ge 52$ & $0,1$ & Every $l\ge 3$ & $l=2$ \\
\bottomrule
\end{tabularx}
\end{center}
Theorem~\ref{thm:global-blowup} realizes  these  endpoints and, with the  known results,  completes the picture for all $n\ge 3$.

We say that $M$ is of \emph{positive Yamabe type} if it admits a smooth metric with positive Yamabe constant, equivalently a metric of positive scalar curvature. The Weyl tensor \emph{vanishes exactly to order $l$ at $p$} if $\nabla^jW(p)=0$ for $j<l$ and $\nabla^lW(p)\ne0$. We call $2$, $4$, and $6$  the endpoint orders in their respective dimension ranges.

\begin{theorem}\label{thm:global-blowup}
Let  $M$ be a closed,  connected,  smooth manifold of dimension $n\ge25$ and of positive Yamabe type,  and let $p\in M$. Let $l\in\N$ satisfy $l\ge 6$ when $25\le n\le 29$, $l \ge 4$ when $30\le n\le 51$, and $l\ge 2$ when $n\ge 52$. Then there exists a single smooth metric $\widehat g$ whose Weyl tensor vanishes exactly to order $l$ at $p$ and satisfies
\[
|\nabla_{\widehat g}^lW_{\widehat g}|>0 \quad\text{everywhere on }M,
\]
such that $[\widehat g]$ contains a sequence of metrics of scalar curvature $4n(n-1)$ whose conformal factors relative to $\widehat g$ blow up at $p$ and whose volumes lie below and converge to that of the round sphere with the same scalar curvature.
\end{theorem}

For orders strictly above the endpoints, the examples of Marques \cite{MarquesBlowup} already give blow-up at $p$ with $|\nabla_{\widehat g}^lW_{\widehat g}(p)|>0$.

\subsection{Proof strategy}
We describe the case $l=2$, $n\ge52$, which exhibits the saddle-point difficulty.
\par\medskip\noindent\textbf{The reduced problem.}
For a single metric $ \widehat g$ to be constructed, we seek solutions near bubbles with centers $x_N+\lambda_N\bar\xi$ and scales $\lambda_N\bar\varepsilon$, where $x_N\to p$ and $\lambda_N\to 0$. Lyapunov--Schmidt reduction solves the equation in the complement of the $n+1$  translation and dilation modes to produce a projected family. Evaluating the Yamabe functional on this family and subtracting its Euclidean bubble value defines the {\it reduced functional}
\[
\widehat F_N(\bar\xi,\bar\varepsilon),
\]
whose critical points give smooth positive solutions of the Yamabe equation for  large $N$. We construct the metric $ \widehat g$ so that the {\it  normalized reduced functionals}  satisfy
\[
\widetilde F_N:=\lambda_N^{-8}\widehat F_N\longrightarrow F \qquad\text{in }C^1
\]
on a fixed neighborhood of $(0,1)$, where the limiting functional $F$ has a nondegenerate saddle point at $(0,1)$ and satisfies $ F(0,1)<0. $

\par\medskip\noindent\textbf{Prescribing the Weyl vanishing order.}
The main geometric difficulty is to construct a single metric $ \widehat g$ with the prescribed second-order Weyl nonvanishing at $p$ and the required limiting functional. In the blow-up constructions of Brendle and Brendle–Marques, the Weyl tensor vanishes to infinite order at the blow-up point. Marques prescribed finite-order Weyl nonvanishing by adding a higher-order perturbation that leaves the limiting functional unchanged. At the endpoint $l=2$, however, it is no longer sufficient for the perturbation responsible for the prescribed Weyl nonvanishing to enter only at higher order in the reduced functional. Instead, in our construction we introduce a perturbation responsible for the prescribed Weyl nonvanishing and arrange that this new perturbation interacts with the localized perturbations already at leading order.

The difference is visible at the level of the polynomial ansatz. In Marques's construction, the tensor governing the limiting functional has the form
\[
H^{\mathrm M}_{ik}(x) =f(|x|^2)\,\Phi(W)_{ik}(x), \qquad \Phi(W)_{ik}(x)=W_{ipkq}x_p x_q,
\]
where $f$ is a polynomial of one variable. Marques then adds a separate higher-order tensor to prescribe the desired finite-order Weyl nonvanishing. Its order is chosen so that it does not affect the limiting functional.

Our endpoint construction enlarges this ansatz by using anisotropic polynomial and two algebraic Weyl tensors. Choose coordinates
\[
x=(t,x'),\qquad t=x_1,\qquad \varrho=|x'|,
\]centered at $p$, and choose nonzero $W,V$ as in
Proposition~\ref{prop:orthogonal-weyl-tensors}. In particular,
\[
W_{1ijk}=V_{1ijk}=0,\qquad \langle W,V\rangle=0,
\]
and their scalar and rank-two quadratic contractions agree. With $\Phi(V)$ defined in the same way as $\Phi(W)$, we take
\[
H=T_{\mathrm{loc}}+T_{\mathrm{glo}},\qquad T_{\mathrm{loc}}=\tau\Phi(W)+\beta t\Phi(V),\qquad T_{\mathrm{glo}}=\varrho^2\Phi(W).
\]
Unlike a radial polynomial, $\varrho^2$ distinguishes the $t$-direction from the transverse variables. The cubic term $\beta t\Phi(V)$ is odd in $t$ and also contributes at leading order. This term will remove the distinguished translation degeneracy and produce a $C^1$-stable saddle point.

The homogeneous tensor $T_{\mathrm{glo}}$ has degree four and prescribes the first nonzero Weyl derivative at the accumulation point. The terms in $T_{\mathrm{loc}}$ have lower degrees and are inserted in small neighborhoods of points $x_N$ approaching the origin.

The localized perturbations are weighted so that their sum is smooth and flat at the origin: every derivative of their sum vanishes there. Consequently, they do not affect any Weyl derivative at $p$, whereas the quartic global term gives
\[
W(p)=0, \qquad \nabla W(p)=0, \qquad \nabla^2W(p)\ne0.
\]
Thus the global perturbation determines the prescribed vanishing order. On the other hand, it is not negligible from the point of view of the blow-up analysis. The points $x_N$ approach the origin along the $t$-axis, while $T_{\mathrm{glo}}$ is independent of $t$. This gives, in the inner region around $x_N$, the exact rescaling
\[
h_0(x_N+y)=\lambda_N^4 H(\lambda_N^{-1}y).
\]
After setting $y=\lambda_N z$, both the global quartic term and the localized terms therefore appear with the same overall factor $\lambda_N^4$. Their quadratic contribution to the reduced functional is of order $\lambda_N^8$.

\par\medskip\noindent\textbf{Removing the axial degeneracy.}
A crucial point is to ensure that the limiting functional has a nondegenerate critical point. The global term $\varrho^2\Phi(W)$ is independent of the distinguished coordinate $t$, so if one used only terms involving $W$, the resulting limiting functional would be invariant under translations in the $t$-direction. Once the remaining parameters are chosen to make $(0,1)$ critical in the transverse translation and dilation variables, this invariance produces a one-dimensional family of degenerate critical points. To break this degeneracy, we introduce the second Weyl tensor $V$ through the odd term
\[
\beta t\Phi(V).
\]
The particular algebraic relations between $W$ and $V$ are chosen so that adding this term breaks the one-dimensional family of degenerate critical points and, after a suitable choice of $\beta$, leaves an isolated nondegenerate critical point. The orthogonality
\[
\langle W,V\rangle=0
\]
eliminates the unwanted first-order interaction in the distinguished translation direction, while the equality of the relevant scalar and rank-two contractions makes the Hessian calculations compatible in the remaining directions. The odd $V$-term produces a strictly negative second derivative in the $t$-translation direction. After choosing $\tau$ appropriately and then taking $\beta>0$ sufficiently small, the remaining translation directions together with the dilation direction form a positive definite block of the Hessian. Consequently,
\[
(\bar\xi,\bar\varepsilon)=(0,1)
\]
is a nondegenerate saddle point of $F$ with exactly one negative direction. The same computation also gives
\[
F(0,1)<0.
\]
The explicit evaluation of the limiting functional and its Hessian is carried out in Section~\ref{sec:auxiliary}, with the required coefficient signs verified in Appendix~\ref{sec:coefficients}. Thus, for $l=2$, the second algebraic Weyl tensor is not merely an auxiliary algebraic device: it is what  removes the axial degeneracy and produces a $C^1$-stable saddle point.

\par\medskip\noindent\textbf{The Euclidean correction and parameter estimates.}
Because the critical point is a saddle rather than a minimum, $C^0$ convergence of the normalized reduced functionals is not sufficient. We need control of their gradients. As in the constructions of Brendle and Marques, we therefore improve the standard bubble by adding a correction obtained by solving the Euclidean linearized equation. The purpose of this correction is to capture the leading effect of the metric perturbation before solving the full nonlinear projected problem. Section~\ref{sec:euclidean-estimates} establishes estimates for this correction, for its parameter derivatives, and for  the corresponding terms in the functional. These estimates identify the limiting functional $F$ and  give the $C^1$ expansion that will later be transferred to the compact manifold.

\par\medskip\noindent\textbf{Global completion.}
We next  place the local metric on an arbitrary closed manifold $M$ of positive Yamabe type. We first choose a background metric with positive Yamabe constant which is Euclidean on a coordinate ball centered at $p$. After making the tensor amplitudes sufficiently small, the local metric constructed above is inserted into this ball without destroying positivity of the Yamabe constant. The local construction already gives
\[
\nabla^2W\ne0
\]
near $p$, but Theorem~\ref{thm:global-blowup} requires the stronger global condition
\[
|\nabla^2W|>0 \qquad\text{everywhere on }M.
\]
Proposition~\ref{prop:weyl-generic} supplies an arbitrarily small perturbation for which $\nabla^2W$ is nonzero everywhere. A cutoff is then used to restore the previously constructed metric near $p$. Since the second Weyl derivative is already nonzero on the intervening annulus, the cutoff can be made without introducing a zero there. This gives a single smooth metric $ \widehat g$ which  agrees with the local metric near $p$ and satisfies
\[
|\nabla_{ \widehat g}^2W_{ \widehat g}|>0 \qquad\text{on }M.
\]

It remains to solve the projected equation for the fixed metric $\widehat g$ on $M$ and to  compare the normalized reduced functional with $F$. A uniform projected inverse permits the nonlinear correction and the required parameter estimates.

Newton's method then gives the nonlinear correction to the approximate solution. Estimates for this correction and, crucially, for its derivatives with respect to the normalized center and scale parameters yield
\[
\lambda_N^{-8}\widehat F_N\longrightarrow F \qquad\text{in }C^1.
\]
A critical point is $C^1$-stable if every sufficiently small $C^1$ perturbation of the functional has a critical point nearby. Here nondegeneracy gives the required stability: since $\nabla F$ has nonzero local Brouwer degree at $(0,1)$, the gradients of the normalized reduced functionals also have a zero nearby for all sufficiently large $N$. Hence $\widehat F_N$ has a critical point
\[
(\bar\xi_N,\bar\varepsilon_N)\longrightarrow(0,1),
\]
and the corresponding corrected bubbles are genuine smooth positive solutions of the Yamabe equation. Their centers tend to $p$ and their scales tend to zero, so they blow up at $p$. Finally, the inequality $F(0,1)<0$ implies that their values of the Yamabe functional, equivalently the volumes of the corresponding constant-scalar-curvature metrics, lie below and converge to the spherical value.

\par\medskip\noindent\textbf{The other endpoints and higher orders.}
This completes the endpoint construction for $l=2$. The endpoint $l=4$ uses the same basic mechanism: a global homogeneous term interacts at leading order with localized perturbations, and a second orthogonal Weyl tensor removes a distinguished translation degeneracy, producing again a nondegenerate saddle point. For $l=6$,  we choose $H$ differently: suitable even $t$-dependent terms already remove the degeneracy, and the limiting functional has a strict local minimum rather than a saddle. Thus the need for the second Weyl tensor and for the saddle-point construction is specific to the endpoint cases $l=2,4$.

For prescribed orders above the endpoint, another modification is needed. In the range $n\ge52$, the endpoint is $l_0=2$. We localize the entire endpoint tensor $H=T_\loc+T_\glo$, multiplying each rescaled copy by the additional factor $\lambda_N^{1/4}$ so that their sum is flat at $p$, and introduce a separate fixed homogeneous tensor of degree $l+2$ which prescribes the desired higher Weyl derivative.  The same $H$ determines $F$, while the new homogeneous term and its interactions contribute only lower-order errors. The leading contribution to the reduced functional is now of order
\[
\lambda_N^{2(l_0+2)+1/2}=\lambda_N^{17/2}.
\]
This separates the two tasks which, at the endpoint, had to be coupled:  the localized copies of $H$ determine $F$, while the separate homogeneous perturbation prescribes the higher Weyl vanishing order. Appendix~\ref{sec:higher-order-proof} carries out this extension.

Sections~\ref{sec:models}--\ref{sec:auxiliary} develop the local metric,  the Euclidean correction and expansion of the functional, and the endpoint critical point. Section~\ref{sec:global} performs the global geometric completion, proves the uniform projected inverse, constructs the projected family, and establishes the $C^1$ convergence of the normalized reduced functionals to the limiting functional. Appendices~\ref{sec:identities}--\ref{sec:KMS} contain the bubble, curvature, integral, and coefficient computations, while Appendix~\ref{sec:higher-order-proof} treats the higher prescribed vanishing orders.

\par\medskip\noindent\textbf{Related results.}
Related results, without attempting an exhaustive list, include the following.  For Branson's fourth-order $Q$-curvature, Wei--Zhao constructed blow-up examples for $n\ge25$ \cite{WeiZhao}, and Li--Xiong proved compactness results in low dimensions \cite{LiXiong}. Assuming positivity of the corresponding Green's function and excluding the round conformal class, Gong--Kim--Wei proved compactness for the fourth-order problem when $5\le n\le24$ and the sixth-order problem when $7\le n\le26$; they also constructed sixth-order blow-up examples for $n\ge27$ \cite{GongKimWei}. For the quotient $Q/R$ of Branson's fourth-order $Q$-curvature by scalar curvature, Ge--Wang--Wei established existence \cite{GeWangWei}. Li--Wang--Wei constructed noncompact families with positive scalar curvature and positive constant $Q/R$ for $n\ge25$; their background metric is smooth for $n\ge26$ and  $C^9$ but not $C^{10}$ for $n=25$ \cite{LiWangWei}. Deng--Lu--Wei constructed smooth blow-up examples for the $\sigma_2$-Yamabe equation for $n\ge27$ \cite{DengLuWei}. Here $\sigma_2$ is the second elementary symmetric function of the eigenvalues, relative to the metric, of the Schouten tensor $\frac{1}{n-2}\bigl(\Ric_g-\frac{R_g}{2(n-1)}g\bigr)$.

\subsection{Notation}
In analytic definitions, $X$ denotes a closed manifold $M$ or $\R^n$, equipped with the indicated metric. Repeated tensor indices are summed from $1$ to $n$ unless stated otherwise; Greek indices in $x'=(x_2,\ldots,x_n)$ range from $2$ to $n$. The prescribed order of the first nonzero Weyl derivative  is $l$. We retain $R_g$, $W_g$, $\nabla_g$, $c(n)$, $\Delta_g$, $L_g$ from (\ref{eq:conformal-law}); $\delta$ is the Euclidean metric, and $\Id$ the identity matrix or operator. On $\R^n$, $\nabla$ denotes Euclidean differentiation.

\begin{itemize}
\item \textbf{Function spaces.} The space
$ \E:=\dot H^1(\R^n) $ is the completion of $C_c^\infty(\R^n)$ under $ \|f\|_\E^2:=\int_{\R^n}|\nabla f|^2\,dx. $ On $M$, we write $H^1(M)=H^1(M,g)$ and use the norm $ \|f\|_{H^1(M,g)}^2:=\int_M(|\nabla_gf|^2+f^2)\,dV_g. $ Their duals are $\E^*$ and $H^{-1}(M)$, respectively; we denote the duality pairing by $\langle G,f\rangle$.

\item \textbf{Functionals.} The quadratic form associated with  $-L_g$  is
$$
Q_g(f):=\int_X\bigl(|\nabla_gf|^2+c(n)R_gf^2\bigr)\,dV_g,
$$
and {\it the Yamabe functional} is $$J_g^X(u):=Q_g(u)-(n-2)^2\int_X|u|^{\frac{2n}{n-2}}\,dV_g.
$$
We omit $X$ when understood. Positive critical points of the Yamabe functional $J_g^M$ are precisely positive solutions of \eqref{eq:yamabe}.
\item \textbf{Equation and errors.} The error term $E_g(v)$,  its linearization at $S$,
and the  nonlinear remainder are
$$
E_g(v):=-L_gv-n(n-2)|v|^{\frac{4}{n-2}}v,\ \ \LL_{g,S}\omega:=-L_g\omega-n(n+2)|S|^{\frac{4}{n-2}}\omega,
$$
\begin{eqnarray*}
\NN_S(\omega)&:=&-[ E_g(S+\omega)-E_g(S)-\LL_{g,S}\omega]\\
&=& n(n-2)\bigl(|S+\omega|^{\frac{4}{n-2}}(S+\omega) -|S|^{\frac{4}{n-2}}S -\frac{n+2}{n-2}|S|^{\frac{4}{n-2}}\omega\bigr).
\end{eqnarray*}
\item \textbf{Algebraic Weyl tensors.} Let $\Wspace$ be the space of algebraic Weyl tensors on $\R^n$,
namely  four-tensors with all the algebraic properties of a Riemannian Weyl tensor. For $W,V\in\Wspace$, set

\begin{equation*}
\langle W,V\rangle:=W_{ijkl}V_{ijkl},\ \ (W\circ V)_{pq}:=W_{pijk}V_{qijk},\ \ \Phi(W)_{ij}(x):=W_{ipjq}x_px_q.
\end{equation*}
Then
\[
|W|^2:= \langle W,W \rangle = \frac{1}{3} \sum_{i,j,k,a}(W_{ikaj}+W_{iakj})^2,
\]
\[
(W\circ V)_{pq} =\frac{1}{3} \sum_{i,j,k}(W_{ikpj}+W_{ipkj})(V_{ikqj}+V_{iqkj}).
\]

\item \textbf{Polynomial tensors.} Let $\Dspace$ be the space of symmetric matrix-valued polynomials $h=(h_{ij})$ on $\R^n$ satisfying
\[
h_{ii}=0,\qquad x_ih_{ij}=0,\qquad \partial_ih_{ij}=0.
\]

The tensor $\Phi(W)$ belongs to $\Dspace$ and is homogeneous of degree two.

\item \textbf{Coordinates.} Write $x=(t,x')$, $t:=x_1$, and $\varrho:=|x'|$ in the polynomial formulas.
We write $B_r(a)$ for a Euclidean ball and $B_r=B_r(0)$, also using this notation for their images in a coordinate chart. The symbol $\mathbf1_A$ denotes the indicator of $A$.
\end{itemize}

\par\medskip\noindent\textbf{AI disclosure.}

Codex and GPT-6 assisted only with the transversality-type Proposition~\ref{prop:weyl-generic} for global Weyl nonvanishing and with manuscript editing and organization. Mathematica was used to compute the lengthy explicit coefficients. All other mathematical content is due entirely to the authors, who independently checked every mathematical statement and proof and take full responsibility for the manuscript.

\par\medskip\noindent\textbf{Acknowledgement.}
Li is partially supported by NSF Grant DMS-2247410 and DMS-2554267. Matos is partially supported by NSF Grant DMS-2247410.

\section{Local metrics with prescribed Weyl derivatives}\label{sec:models}
We construct endpoint metrics with compactly supported perturbations of the Euclidean metric. The metric will later be inserted into a coordinate ball of $M$.

\subsection{Algebraic Weyl tensors and homogeneous polynomials}
Recall that $\Phi(W)_{ij}(x)=W_{ipjq}x_px_q$ and $(W\circ V)_{pq}=W_{pijk}V_{qijk}$. The following pair has matching quadratic contractions but is orthogonal.
\begin{proposition}[{Orthogonal Weyl tensors with matching contractions}]\label{prop:orthogonal-weyl-tensors}
For every $n\ge5$ there are nonzero $W,V\in\Wspace$ such that
\begin{equation}\label{eq:pair-structural}
W_{1ijk}=V_{1ijk}=0, \ \  \ip WV=0, \ \ |W|^2=|V|^2,\ \ W\circ W=V\circ V.
\end{equation}

Moreover, if $e_1,\ldots,e_n$ denote the standard orthonormal basis of $\mathbb R^n$ and $ E_4:=\operatorname{span}\{e_2,e_3,e_4,e_5\}, $ then $W$ and $V$ may be chosen so that
\begin{equation}\label{eq:pair-normalization}
W\circ W=V\circ V= 6\sum_{\alpha=2}^5 e_\alpha\otimes e_\alpha.
\end{equation}

\end{proposition}
\begin{proof}

Orient $E_4$ and split $ \Lambda^2E_4=\Lambda^2_+E_4\oplus\Lambda^2_-E_4 $ into the $+1$ and $-1$ eigenspaces of the Hodge star. Write $e_{ij}:=e_i\wedge e_j$.

In the orthonormal bases $ \omega_1^\pm= (e_{23}\pm e_{45})/\sqrt2$, $\omega_2^\pm=(e_{24}\mp e_{35})/\sqrt2$, and $\omega_3^\pm=(e_{25}\pm e_{34})/\sqrt2, $ define $ W=2\,\omega_1^+\otimes\omega_1^+-\omega_2^+\otimes\omega_2^+-\omega_3^+\otimes\omega_3^+$ and $ V=2\,\omega_1^-\otimes\omega_1^--\omega_2^-\otimes\omega_2^--\omega_3^-\otimes\omega_3^-. $ Extend them by zero on $E_4^\perp$.

The contractions can be computed from the two-form representation. Writing $(a_1,a_2,a_3)=(2,-1,-1)$, their components are $ W_{ijkl}=\sum_{r=1}^3 a_r(\omega_r^+)_{ij}(\omega_r^+)_{kl}$ and $V_{ijkl}=\sum_{r=1}^3 a_r(\omega_r^-)_{ij}(\omega_r^-)_{kl}. $

Properties  \eqref{eq:pair-normalization} and  \eqref{eq:pair-structural} follow immediately.
\end{proof}

A common rescaling of $W,V$ preserves \eqref{eq:pair-structural}. We choose its size after fixing the polynomial coefficients in Section~\ref{sec:auxiliary}.

In this section $l\in\{2,4,6\}$ and $n\ge25$. Fix the pair $W$ and $V$ as above and let
\begin{equation*}
T_\glo(x):=\varrho^l\Phi(W)(x).
\end{equation*}
It has degree $l+2$ and is independent of $t$. It will prescribe the Weyl vanishing order at the origin. The lower-degree tensor $T_\loc$ below permits adjustment of the reduced functional. For integers $a,b\ge 0$, write
\begin{equation*}
H^{a,b}_{ij}:=t^{2a}\varrho^{2b}\Phi(W)_{ij}.
\end{equation*}
Set
\begin{align*}
H^{(2m+2)}&:=\sum_{a+b=m}\tau_{a,b}H^{a,b}, \qquad 0\le m\le l/2-1,\\
H^{(3)}&:=\beta t\Phi(V),\qquad T_\loc:=\sum_{m=0}^{l/2-1}H^{(2m+2)}+H^{(3)},\qquad H:=T_\loc+T_\glo.
\end{align*}
The constants $\tau_{a,b}$ and $\beta$ are chosen in Section~\ref{sec:auxiliary}. For homogeneous expansions below, set $H^{(l+2)}=T_\glo$ and set every unlisted homogeneous component $H^{(d)}$ to zero. Each tensor belongs to $\Dspace$. Indeed, $\Phi(W)$ is trace-free, divergence-free and satisfies $x_i\Phi(W)_{ij}=0$. The identities $\Phi(W)_{1j}=0$ and $\sum_{\alpha=2}^n x_\alpha\Phi(W)_{\alpha j}=0$ show, by the product rule, that $t^{2a}\varrho^{2b}\Phi(W)$ is also divergence-free. The same reasoning applies to the $V$ term.

\subsection{Metric construction}
Choose $\eta\in C_c^\infty((-1/2,1/2))$, nonnegative and equal to one on $[-1/4,1/4]$. For
\begin{equation*}
x_N:=(N^{-1},0),\qquad \lambda_N:=2^{-N},\qquad \rho_N:=(4N^2)^{-1},
\end{equation*}
define $g_0:=\exp(h_0)$, where
\begin{equation}\label{eq:h0}
h_0(x):=\eta(|x|^2)T_\glo(x) +\sum_{N=N_0}^{\infty}\eta(N^4|x-x_N|^2) \lambda_N^{l+2}T_\loc(\lambda_N^{-1}(x-x_N)).
\end{equation}
For large $N_0$, all supports are contained in $B_1$.
\begin{lemma}[{Smoothness and exact local rescaling}]\label{lem:smooth-bump-series}
The infinite series in \eqref{eq:h0} converges in $C^\infty$, and its sum vanishes to infinite order at $0$. For all sufficiently large $N$,
\[
h_0(x_N+y)=\lambda_N^{l+2}H(\lambda_N^{-1}y), \qquad |y|<\rho_N.
\]
Every fixed $C^k$ norm of $h_0$ can be made arbitrarily small by a common rescaling of $W,V$.
\end{lemma}

On $B_{\rho_N}(x_N)$ the $N$th cutoff and the global cutoff equal one, and every other local cutoff is zero. Independence of $T_\glo$ from $t$, together with its homogeneity, gives the exact rescaling above. Since $H\in\Dspace$, it also gives, with $\lambda=\lambda_N$ and $\rho=\rho_N$,
\[
y_i(h_0)_{ij}(x_N+y)=0,\qquad \partial_{y_i}(h_0)_{ij}(x_N+y)=0\quad\text{on }B_\rho.
\]

The common rescaling of the two Weyl tensors makes $\|h_0\|_{C^2}$ sufficiently small. For all sufficiently large $N$, set
\begin{equation}\label{eq:translated}
\bar h(y):=h_0(x_N+y),\qquad \bar g:=\exp(\bar h),\qquad H^\lambda(y):=\lambda^{l+2}H(\lambda^{-1}y).
\end{equation}
Thus, $\bar h=H^\lambda$ on $B_\rho$. The original cutoff supports lie in a ball of radius strictly smaller than one, so $\bar h$ is also supported in $B_1$ for large $N$.

\begin{lemma}[{Weyl vanishing order}]\label{lem:weyl-vanishing-order}
At the accumulation point, $g_0$ satisfies
\begin{equation*}
\nabla_{g_0}^jW_{g_0}(0)=0\quad(0\le j<l),\qquad \nabla_{g_0}^lW_{g_0}(0)\ne0.
\end{equation*}
More precisely, with the pair supported on $E_4$, $\nabla_{e_6}^{\,l}W_{g_0}(0)=-3l!W$.
\end{lemma}
\begin{proof}

Lemma~\ref{lem:smooth-bump-series} gives $h_0-T_\glo$ flat at $0$, and hence $h_0=O(|x|^{l+2})$. The linearized Riemann tensor in our convention is
\begin{equation}\label{eq:linear-curvature}
(D\Rm_\delta[h])_{ijka} =\frac{1}{2}(\partial_{jk}h_{ia}+\partial_{ia}h_{jk} -\partial_{ja}h_{ik}-\partial_{ik}h_{ja}).
\end{equation}
Since $\partial_{ab}\Phi(W)_{ij}=W_{iajb}+W_{ibja}$, the curvature symmetries and first Bianchi identity give
\begin{equation}\label{eq:linear-Weyl}
D\Rm_\delta[\Phi(W)]=-3W,\qquad D\Ric_\delta[\Phi(W)]=DR_\delta[\Phi(W)]=0.
\end{equation}
Thus $DW_\delta[\Phi(W)]=-3W$. Since $l$ is even,
\[
\partial_6^lT_\glo=\partial_6^l(\varrho^l\Phi(W))=l!\Phi(W).
\]
The nonlinear curvature terms are $O(|h_0||\partial^2h_0|+|\partial h_0|^2)=O(|x|^{2l+2})$, so their derivatives through order $l$ vanish at $0$. The Christoffel symbols satisfy $\Gamma_{g_0}=O(|x|^{l+1})$, and coordinate and covariant curvature derivatives agree through order $l$ at the origin. Constant-coefficient differentiation of \eqref{eq:linear-Weyl} therefore gives $\nabla_{e_6}^{\,l}W_{g_0}(0)=-3l!W\ne0$. All lower derivatives vanish by homogeneity. The same assertions hold for any smooth metric agreeing with $g_0$ near $0$.

\end{proof}

\section{ Euclidean correction and functional estimates}\label{sec:euclidean-estimates}
We construct the Euclidean linear correction that captures the leading effect of the metric perturbation. We estimate this correction and its first parameter derivatives, and identify the limiting functional that will govern the reduction on $M$.

Fix the metric $g_0=\exp(h_0)$  in \eqref{eq:h0}. Throughout this section, $l\in\{2,4,6\}$ and $n\ge25$. For the analytic estimates in this section, Section~\ref{sec:global}, and Appendix~\ref{sec:higher-order-proof}, fix
\[
m:=\frac{n-2}{2},\qquad p:=\frac{n+2}{n-2},\qquad q:=\frac{2n}{n+2}.
\]
Let $\lambda=2^{-N}$, $\rho=(4N^2)^{-1}$, and use the translated metric \eqref{eq:translated}.

\subsection{The linear correction and error estimates}\label{subsec:linear-correction}
For $(\xi,\varepsilon)\in\R^n\times(0,\infty)$, let
\begin{equation}\label{TheBubble}
u_{(\xi,\varepsilon)}(x):= \left(\frac{\varepsilon}{\varepsilon^2+|x-\xi|^2}\right)^{\frac{n-2}{2}}.
\end{equation}
Direct computations give
\[
-\Delta u_{(\xi,\varepsilon)} =n(n-2)u_{(\xi,\varepsilon)}^{\frac{n+2}{n-2}},
\]
and
\begin{equation}
J_\sph :=J_\delta^{\R^n}(u_{(\xi,\varepsilon)})\equiv 2(n-2) \int_{\R^n}u_{(\xi,\varepsilon)}^{\frac{2n}{n-2}}\dd \equiv 2(n-2) \int_{\R^n}u_{(0, 1)}^{\frac{2n}{n-2}}\dd.
\label{eq:spherical-level}
\end{equation}

Let $Z_{(\xi,\varepsilon,0)}:=\partial_\varepsilon u_{(\xi,\varepsilon)}$ and $Z_{(\xi,\varepsilon,j)}:=\partial_{\xi_j}u_{(\xi,\varepsilon)}$, $1\le j\le n$. These parameter derivatives are called the Jacobi fields of the bubble. They satisfy
\begin{equation*}
-\Delta Z_{(\xi,\varepsilon,j)} -n(n+2)u_{(\xi,\varepsilon)}^{\frac{4}{n-2}}Z_{(\xi,\varepsilon,j)}=0.
\end{equation*}

The bubble pairings and norm bounds used below are recorded in Appendix~\ref{subsec:bubble-identities}. Set
\[
\varphi_{(\xi,\varepsilon,0)} :=-\frac{2\varepsilon}{n-2} u_{(\xi,\varepsilon)}^{\frac{4}{n-2}} Z_{(\xi,\varepsilon,0)},\ \ \ \varphi_{(\xi,\varepsilon,j)} :=\frac{2\varepsilon}{n-2} u_{(\xi,\varepsilon)}^{\frac{4}{n-2}} Z_{(\xi,\varepsilon,j)}, \ \ 1\le j\le n.
\]

Define the closed subspace
\begin{equation*}
\E_{(\xi,\varepsilon)}:= \left\{f\in\E:\int_{\R^n}\varphi_{(\xi,\varepsilon,j)}f\dd=0, \quad0\le j\le n\right\}.
\end{equation*}
We have that $u_{(\xi,\varepsilon)}\in\E_{(\xi,\varepsilon)}.$

Define $\LL^\perp_{\delta,(\xi,\varepsilon)}: \E_{(\xi,\varepsilon)}\to\E_{(\xi,\varepsilon)}^*$ by setting, for $f,\psi\in\E_{(\xi,\varepsilon)}$,
\begin{equation*}
\left\langle \LL^\perp_{\delta,(\xi,\varepsilon)}f,\psi \right\rangle := \int_{\R^n}\bigl( \partial_if\partial_i\psi -n(n+2)u_{(\xi,\varepsilon)}^{\frac{4}{n-2}}f\psi \bigr)\dd.
\end{equation*}

\par\medskip\noindent\textbf{Euclidean projected inverse.}
It is well known, see e.g. \cite[proof of Proposition~4]{Brendle},   that the operator $\LL^\perp_{\delta,(\xi,\varepsilon)}$ is an isomorphism satisfying
\begin{equation}
\label{flat-projected-inverse}
\|f\|_\E \le C(n) \|\LL^\perp_{\delta,(\xi,\varepsilon)} f\|_{\E_{(\xi,\varepsilon)}^*},\quad \text{for all}\ f\in \E_{(\xi, \varepsilon)}.
\end{equation}

The identities $H_{ii}=0$, $x_iH_{ij}=0$, and $\partial_iH_{ij}=0$ imply the following orthogonality.
\begin{lemma}[{Orthogonality of the leading source}]\label{lem:source-orthogonality}
For every $H\in\Dspace$ with $\deg H\le8$ and every Jacobi field,
\begin{equation*}
\int_{\R^n}Z_{(\xi,\varepsilon,j)}H_{ab} \partial_{ab}u_{(\xi,\varepsilon)}\dd=0, \qquad0\le j\le n.
\end{equation*}
\end{lemma}
\begin{proof}
Write $u:=u_{(\xi,\varepsilon)}$. Direct differentiation yields
\begin{equation}
\label{eq:stress-identity}
\partial_iu\partial_ju-c(n)\partial_{ij}(u^2) = (n-2)^2 u^2 [2(n-1)(\varepsilon^2+|x-\xi|^2)]^{-1}
\delta_{ij}.
\end{equation}
Contract with $H$ and integrate twice by parts. Its trace and divergence vanish, so
\[
\int_{\R^n}H_{ij}\partial_iu\partial_ju\dd =c(n)\int_{\R^n}\partial_{ij}H_{ij}u^2\dd=0.
\]

Since the degree is at most eight, the boundary terms at radius $R$ are $O(R^{10-n})$, tending to zero for $n\ge25$. The same decay permits differentiation in every bubble parameter. Integrating once more by parts yields the lemma.

\end{proof}
Let $H$ be the polynomial tensor of Section~\ref{sec:models}, and recall that $H^\lambda(x)=\lambda^{l+2}H(\lambda^{-1}x)$. For every $(\xi,\varepsilon)\in\R^n\times(0,\infty)$, the function $H_{ij}\partial_{ij}u_{(\xi,\varepsilon)}$ belongs to $L^{\frac{2n}{n+2}}(\R^n)$, since $\deg H=l+2$ and $n>2l+6$. Define $\bar w$ by
\begin{equation}\label{eq:w-equation}
\left\langle\LL^\perp_{\delta,(\xi,\varepsilon)} \bar w_{(\xi,\varepsilon)},\psi\right\rangle =-\int_{\R^n}H_{ij}\partial_{ij}u_{(\xi,\varepsilon)}\psi\dd, \qquad \psi\in\E_{(\xi,\varepsilon)}.
\end{equation}\begin{lemma}[{The Euclidean linear correction}]\label{lem:corrections-full-equations}
For every $(\xi,\varepsilon)\in\R^n\times(0,\infty)$, equation~\eqref{eq:w-equation} has a unique solution $\bar w_{(\xi,\varepsilon)}\in\E_{(\xi,\varepsilon)}$. This solution is smooth and satisfies the full linearized equation
\begin{equation}\label{eq:unit-full}
\Delta\bar w_{(\xi,\varepsilon)} +n(n+2)u_{(\xi,\varepsilon)}^{\frac{4}{n-2}} \bar w_{(\xi,\varepsilon)} =H_{ij}\partial_{ij}u_{(\xi,\varepsilon)}, \qquad \text{in}\ \R^n.
\end{equation}
\end{lemma}
\begin{proof}
Existence and uniqueness follow from the projected inverse \eqref{flat-projected-inverse}. We suppress the bubble parameters $(\xi,\varepsilon)$ in the rest of the proof. Equation~\eqref{eq:w-equation} implies that there are constants $b_0,\ldots,b_n$ such that
\[
\Delta\bar w+n(n+2)u^{\frac{4}{n-2}}\bar w =H_{ij}\partial_{ij}u+\sum_{j=0}^n b_j\varphi_j, \qquad \text{in}\ \R^n.
\]
Testing the above  with $Z_k\in\E$, and using $\Delta Z_k+n(n+2)u^{\frac{4}{n-2}}Z_k=0$ together with Lemma~\ref{lem:source-orthogonality} applied to $H$, we obtain
\[
\begin{aligned}
\sum_{j=0}^n b_j\int_{\R^n}\varphi_jZ_k\dd &=\int_{\R^n}\bar w\bigl(\Delta Z_k+ n(n+2)u^{\frac{4}{n-2}}Z_k\bigr)\dd\\
&\quad-\int_{\R^n}H_{ij}\partial_{ij}u\,Z_k\dd =0,\qquad 0\le k\le n.
\end{aligned}
\]
By \eqref{eq:dual-Jacobi} in Lemma~\ref{lem:bubble-pairings}, the pairing matrix is diagonal with nonzero entries. Hence $b_0=\cdots=b_n=0$, proving \eqref{eq:unit-full} weakly. Local elliptic regularity makes $\bar w$ smooth, so the equation holds pointwise.
\end{proof}

\par\medskip\noindent\textbf{Normalized parameters}
For the estimates below, set
\[
\Omega:=B_1(0)\times(1/2,2),\qquad U:=B_2(0)\times(1/4,4).
\]Recall that $ \lambda = 2^{-N} $and write,
\[
(\xi,\varepsilon)=\lambda(\bar\xi,\bar\varepsilon),\qquad (\bar\xi,\bar\varepsilon)\in U.
\]
All parameter-dependent families are defined and differentiated on $U$, and their uniform estimates are taken over $\overline\Omega$. We use $\theta$ for any of $\xi_1,\ldots,\xi_n,\varepsilon$ and $\bar\theta$ for its normalized counterpart, so $\partial_{\bar\theta}=\lambda\partial_\theta$. 

For a Banach space $X$, the notation $C^1(U,X)$ refers to continuously differentiable maps on the normalized parameter domain $U$. For their restrictions to $\overline\Omega$, we use the norm
\[
\|f\|_{C^1(\overline\Omega,X)} :=\sup_{(\bar\xi,\bar\varepsilon)\in\overline\Omega} \bigl(\|f\|_X+\sum_{j=1}^n\|\partial_{\bar\xi_j}f\|_X +\|\partial_{\bar\varepsilon}f\|_X\bigr).
\]
All bounds in this norm are uniform in $N$. Remainder estimates below are stated as bounds or limits for these norms.

\par\medskip\noindent\textbf{The scaled linear correction.}
Define
\begin{equation}\label{eq:exact-scaling}
w_{(\xi,\varepsilon)}(y) :=\lambda^{l+2-\frac{n-2}{2}} \bar w_{(\bar\xi,\bar\varepsilon)}(\lambda^{-1}y), \qquad y\in\R^n.
\end{equation}
By a change of variables, $w_{(\xi,\varepsilon)}\in\E_{(\xi,\varepsilon)}$ satisfies \eqref{eq:w-equation} with $\bar w$ replaced by $w$ and $H$ by $H^\lambda$. By  \eqref{eq:unit-full},
\begin{equation}\label{eq:w-full}
\Delta w_{(\xi,\varepsilon)} +n(n+2)u_{(\xi,\varepsilon)}^{\frac{4}{n-2}} w_{(\xi,\varepsilon)} =H^\lambda_{ij}\partial_{ij}u_{(\xi,\varepsilon)},\qquad \text{in}\ \ \R^n.
\end{equation}

The correction $w_{(\xi,\varepsilon)}$ is chosen to cancel the leading source $H^\lambda_{ij}\partial_{ij}u_{(\xi,\varepsilon)}$ in the bubble error $E_{\bar g}(u_{(\xi,\varepsilon)})$. We next estimate this error and the remainder after subtraction.
\begin{proposition}[{Bubble error and leading source}]\label{prop:bubble-estimates}
There exists a constant $C>0$, such that for all $(\bar\xi,\bar\varepsilon)\in\overline\Omega$ and $s=0,1$,
\begin{align}
\|\partial_\theta^s E_{\bar g}(u_{(\xi,\varepsilon)})\|_{L^q(\R^n)} &\le C \lambda^{-s}\bigl(\lambda^{l+2}+(\lambda/\rho)^m\bigr),\label{eq:res-bubble}\\
\|\partial_\theta^s[E_{\bar g}(u_{(\xi,\varepsilon)})-H^\lambda_{ij}\partial_{ij}u_{(\xi,\varepsilon)}]\|_{L^q(\R^n)} &\le C\lambda^{-s}\bigl(\lambda^{p(l+2)}+(\lambda/\rho)^m\bigr).\label{eq:res-improved}
\end{align}
\end{proposition}
\begin{proof}

For $A>0$ and $0<\lambda\le\rho<1$, we have
\begin{equation}\label{eq:radial-master}
\|\lambda^{\frac{n-2}{2}}(\lambda+|x|)^{A-n}\|_{L^{\frac{2n}{n+2}}(B_\rho)}
\le C\begin{cases}
\lambda^A,&n>2A+2,\\
\lambda^A(1+\log(\frac{\rho}{\lambda}))^{\frac{n+2}{2n}},&n=2A+2,\\
\lambda^{\frac{n-2}{2}}\rho^{A-\frac{n-2}{2}},&n<2A+2.
\end{cases}
\end{equation}

Write $u:=u_{(\xi,\varepsilon)}$ in this proof.

For $0\le a\le3$ and $s=0,1$, the explicit bubble satisfies
\begin{equation}\label{eq:bubble-pointwise}
|\partial_\theta^s\nabla^au_{(\xi,\varepsilon)}(x)| \le C\lambda^{m-s}(\lambda+|x|)^{2-n-a}.
\end{equation}
On $B_\rho$, $\bar h=H^\lambda=\sum_{d\le (l+2)}\lambda^{l+2-d}H^{(d)}$ is divergence-free and satisfies
\[
|\nabla^a\bar h(x)|\le C(\lambda+|x|)^{l+2-a},\qquad0\le a\le2.
\]
The scalar-curvature expansion in Lemma~\ref{lem:scalar-expansion} gives
\begin{equation}\label{eq:core-R}
|R_{\bar g}|\le C(|\bar h|^2|\partial^2\bar h|+|\partial\bar h|^2) \le C(\lambda+|x|)^{2(l+2)-2}.
\end{equation}
Since $\tr\bar h=0$, we have $\det\bar g=1$, and thus  obtain
\begin{align*}
-E_{\bar g}(u)&=\partial_i((\bar g^{ij}-\delta^{ij})\partial_ju)-c(n)R_{\bar g}u,\\
-E_{\bar g}(u)+H^\lambda_{ij}\partial_{ij}u &=\partial_i((\bar g^{ij}-\delta^{ij}+H^\lambda_{ij})\partial_ju)-c(n)R_{\bar g}u.
\end{align*}
The second equation uses $\partial_iH^\lambda_{ij}=0$. Matrix exponentiation gives $|\bar g^{-1}-\Id+\bar h|\le C|\bar h|^2$ and $|\partial(\bar g^{-1}-\Id+\bar h)|\le C|\bar h||\partial\bar h|$. Consequently on $B_\rho$,
$$
|E_{\bar g}(u)|\le C\lambda^m(\lambda+|x|)^{l+2-n},\ \ \ |E_{\bar g}(u)-H^\lambda_{ij}\partial_{ij}u| \le C\lambda^m(\lambda+|x|)^{2(l+2)-n}.
$$
It follows, using ~\eqref{eq:radial-master}, that $ \|E_{\bar g}(u_{(\xi,\varepsilon)})\|_{L^{2n/(n+2)}(B_\rho)} \le C\lambda^{l+2}$, and $ \|E_{\bar g}(u_{(\xi,\varepsilon)})-H^\lambda_{ij}\partial_{ij}u_{(\xi,\varepsilon)}\|_{L^{2n/(n+2)}(B_\rho)} \le C\lambda^{p(l+2)}$.

Using (\ref{eq:bubble-pointwise}), we have $\|E_{\bar g}(u_{(\xi,\varepsilon)})\|_{  L^q(\rho\le |x|\le 1) } \le C(\lambda/\rho)^m. $ For large $N$, $\bar h=0$ on $|x|\ge1$, so $E_{\bar g}(u_{(\xi,\varepsilon)})=0$ there.

The polynomial source $H^\lambda_{ij}\partial_{ij}u_{(\xi,\varepsilon)}$ satisfies, for each degree $d\le (l+2)<m$,
\begin{equation*}
\|\lambda^{l+2-d}H^{(d)}_{ij}\partial_{ij}u\|_{L^q(\R^n\setminus B_\rho)} \le C\lambda^{m+(l+2)-d}\rho^{d-m} \le C(\lambda/\rho)^m.
\end{equation*}

This proves \eqref{eq:res-bubble}--\eqref{eq:res-improved} for $s=0$. Applying $\partial_\theta $ to the expression of $E_{\bar g}(u)$ and $E_{\bar g}(u)-H^\lambda_{ij}\partial_{ij}u$, we obtain, using (\ref{eq:bubble-pointwise}), the estimates for $s=1$.
\end{proof}

\subsection{Estimates for the linear correction}

The following lemma provides the required $C^1$ parameter dependence and derivative estimates for solutions of the projected linear equation. The projection onto $\E_{(\xi,\varepsilon)}$ is
\begin{equation}\label{eq:projection}
T_{(\xi,\varepsilon)}z :=z-\sum_{j=0}^n \frac{\int_{\R^n}\varphi_{(\xi,\varepsilon,j)}z\dd}{\|\varphi_{(\xi,\varepsilon,j)}\|_2^2} \varphi_{(\xi,\varepsilon,j)}.
\end{equation}
Equation~\eqref{eq:dual-dual} shows that these projections identify nearby subspaces.

\begin{lemma}[{Parameter dependence of the projected inverse}]\label{lem:parameter-derivative-inverse}
Let $(\xi,\varepsilon)\mapsto f_{(\xi,\varepsilon)}$ be a $C^1$ map from $\R^n\times(0,\infty)$ into $L^{\frac{2n}{n+2}}(\R^n)$. Then the family $\psi_{(\xi,\varepsilon)} :=(\LL^\perp_{\delta,(\xi,\varepsilon)})^{-1} f_{(\xi,\varepsilon)}$ is $C^1$ with values in $\E$. For $\theta\in\{\xi_1,\ldots,\xi_n,\varepsilon\}$,
\begin{equation*}
\|\partial_\theta\psi_{(\xi,\varepsilon)}\|_\E \le C\bigl( \|\partial_\theta f_{(\xi,\varepsilon)}\|_{L^{\frac{2n}{n+2}}} +\varepsilon^{-1} \|f_{(\xi,\varepsilon)}\|_{L^{\frac{2n}{n+2}}} \bigr).
\end{equation*}
\end{lemma}

\begin{proof}
Write $a:=(\xi,\varepsilon)$ and fix $a_0$. The space $\E_a$ depends on $a$. We first transfer the equation to the fixed space $\E_{a_0}$. By \eqref{eq:dual-dual}, $T_a$ in \eqref{eq:projection} is a projection onto $\E_a$. The explicit fields and \eqref{eq:field-norms} show that $a\mapsto T_a$ is $C^1$ in the operator norm on $\E$, with $\|T_a\|\le C$ and $\|\partial_\theta T_a\|\le C\varepsilon^{-1}$. For $a$ sufficiently close to $a_0$, the operator $\widetilde T_a:=T_a|_{\E_{a_0}}:\E_{a_0}\to\E_a$ satisfies
\begin{equation}\label{eq:fixed-space-U}
\|\widetilde T_a\|+\|\widetilde T_a^{-1}\|\le C,\qquad \|\partial_\theta \widetilde T_a\|\le C\varepsilon^{-1}.
\end{equation}

On $\E_{a_0}$ define
\[
\langle A_az,\eta\rangle :=\left\langle \LL^\perp_{\delta,a}\widetilde T_az,\widetilde T_a\eta \right\rangle, \qquad \langle F_a,\eta\rangle :=\int_{\R^n}f_a\widetilde T_a\eta\dd.
\]
Estimate~\eqref{flat-projected-inverse} bounds $A_a^{-1}$. Since $|\partial_\theta u_a^{\frac{4}{n-2}}| \le C\varepsilon^{-1}u_a^{\frac{4}{n-2}}$, H\"older's inequality and \eqref{eq:fixed-space-U} imply
\[
\|\partial_\theta A_a\|\le C\varepsilon^{-1},\qquad \|\partial_\theta F_a\| \le C\bigl( \|\partial_\theta f_a\|_{\frac{2n}{n+2}} +\varepsilon^{-1}\|f_a\|_{\frac{2n}{n+2}} \bigr).
\]
The identity $\partial_\theta A_a^{-1} =-A_a^{-1}(\partial_\theta A_a)A_a^{-1}$ and the product rule for $\psi_a=\widetilde T_aA_a^{-1}F_a$ prove the assertion. Translation and dilation make the bounds uniform at the scale of $\varepsilon$.
\end{proof}

We now apply the parameter-dependence lemma to the scaled linear correction $w_{(\xi,\varepsilon)}$, obtaining spatial decay and bounds for its first parameter derivatives.

\begin{proposition}[{Estimates for $w$}]\label{prop:w-decay}
Let $w=w_{(\xi,\varepsilon)}$ be defined by \eqref{eq:exact-scaling}, where $(\xi,\varepsilon)=\lambda(\bar\xi,\bar\varepsilon)$. There is a constant $C>0$, independent of $N$, such that, for all $(\bar\xi,\bar\varepsilon)\in\overline\Omega$, $a\in\{0,1,2\}$, and $\theta\in\{\xi_1,\ldots,\xi_n,\varepsilon\}$,
\begin{align}
|\nabla^aw(x)|+\lambda|\partial_\theta\nabla^aw(x)| &\le C\lambda^{\frac{n-2}{2}}(\lambda+|x|)^{l+4-n-a}, \quad x\in \R^n,\label{eq:w-decay}\\
\|w\|_\E&\le C\lambda^{l+2}.
\label{eq:w-energy}
\end{align}
\end{proposition}
\begin{proof}
Write $u=u_{(\xi,\varepsilon)}$. Recall that $w=-(\LL^\perp_{\delta,(\xi,\varepsilon)})^{-1} (H^\lambda_{ij}\partial_{ij}u)$ and satisfies \eqref{eq:w-full}. The expression of $u$ gives
\begin{equation}\label{AAA}
|\nabla^a(H^\lambda_{ij}\partial_{ij}u)| \le C\lambda^{\frac{n-2}{2}} (\lambda+|x|)^{l+2-n-a},\qquad a=0,1,2.
\end{equation}
Since $n>2l+6$,  we have $\|H^\lambda_{ij}\partial_{ij}u\|_{L^{2n/(n+2)}(\R^n)} \le C\lambda^{l+2}. $ Thus, by  (\ref{flat-projected-inverse}), $ \|w\|_\E \le C\|H^\lambda_{ij}\partial_{ij}u\|_{L^{2n/(n+2)}(\R^n)}. $ Estimate \eqref{eq:w-energy} is established.

The estimate for $|\nabla^a w(x)|$ in \eqref{eq:w-decay} follows by the argument of \cite[Proposition~9]{Brendle}, using (\ref{AAA}). Indeed, after rescaling to $\lambda=1$,  estimate~\eqref{eq:w-energy} and interior estimates of $w$ give its decay of order $(1+|x|)^{-(n-2)/2}$. The weighted iteration in that proof improves the decay to $(1+|x|)^{l+4-n}$, since $(n-2)/2<n-l-4<n-2$. This establishes the estimate for $|w(x)|$ in (\ref{eq:w-decay}). Standard derivative estimates then give the estimate for $|\nabla^aw(x)|$.

It remains to estimate the parameter-derivative term in \eqref{eq:w-decay}. By Lemma~\ref{lem:parameter-derivative-inverse}, we have $\|\partial_\theta w\|_\E\le C\lambda^{l+1}$. Differentiating \eqref{eq:w-full} gives
\[
\Delta(\partial_\theta w) +n(n+2)u^{\frac4{n-2}}\partial_\theta w ={}H^\lambda_{ij}\partial_{ij}\partial_\theta u -n(n+2)\partial_\theta(u^{\frac4{n-2}})w.
\]
By \eqref{eq:w-decay}, the right-hand side of the above and its spatial derivatives of order $a=0,1, 2$ are bounded by $ C\lambda^{\frac{n-4}{2}} (\lambda+|x|)^{l+2-n-a}. $ Thus the same arguments as above, applied to $\partial_\theta w$, finishes the proof of \eqref{eq:w-decay}.
\end{proof}
To use the linear correction in the metric $\bar g$, we must also estimate the error caused by replacing $\Delta$ with $L_{\bar g}$.

\begin{lemma}[{Error of the linear correction in the perturbed metric}]\label{lem:linear-correction-error}
There exists a constant $C>0$ such that, for all $(\bar\xi,\bar\varepsilon)\in\overline\Omega$ and $s=0,1$, the correction $w_{(\xi,\varepsilon)}$ defined by \eqref{eq:exact-scaling} satisfies
\begin{equation}\label{eq:linear-correction-error}
\|\partial_\theta^s[(L_{\bar g}-\Delta)w_{(\xi,\varepsilon)}]\|_{L^{\frac{2n}{n+2}}} \le C\lambda^{-s}\bigl( \lambda^{\frac{(n+2)(l+2)}{n-2}} +(\lambda/\rho)^{\frac{n-2}{2}}\bigr).
\end{equation}
\end{lemma}
\begin{proof}
Write $w:=w_{(\xi,\varepsilon)}$, and set
\[
B_1:=(L_{\bar g}-\Delta)w =\partial_i((\bar g^{ij}-\delta^{ij})\partial_jw)-c(n)R_{\bar g}w.
\]
On $B_\rho$, \eqref{eq:core-R} and Proposition~\ref{prop:w-decay} give $ |B_1|\le C\lambda^m(\lambda+|x|)^{2(l+2)-n}. $

Applying \eqref{eq:radial-master} with $A=2(l+2)$ gives $\|B_1\|_{L^q(B_\rho)}\le C\lambda^{p(l+2)}$, since $2(l+2)>p(l+2)$ and $m\ge p(l+2)$. On $\rho\le|x|\le1$, \eqref{eq:w-decay} gives
\[
|B_1|\le C\lambda^m|x|^{l+2-n},\qquad \|B_1\|_{L^q(B_1(0)\setminus B_\rho)}\le C\lambda^m\rho^{l+2-m} \le C(\lambda/\rho)^m.
\]
Since $B_{1}(x)$ vanishes outside $B_1(0)$, this proves \eqref{eq:linear-correction-error} for $s=0$. The same argument, using the parameter-derivative bound in \eqref{eq:w-decay},  gives \eqref{eq:linear-correction-error} for $s=1$.
\end{proof}

\subsection{The geometric terms in the limiting functional}
All integrals in this subsection are over $\R^n$ with measure $dx$, unless a domain is displayed. Define the limiting functional by
\begin{equation}\label{eq:auxiliary-function}
\begin{split}
F(\xi,\varepsilon):={}&\frac{1}{2}\int_{\R^n} H_{ik}H_{jk}\partial_iu_{(\xi,\varepsilon)} \partial_ju_{(\xi,\varepsilon)}\dd\\
&-\frac{c(n)}{4}\int_{\R^n}(\partial_kH_{ij})^2 u_{(\xi,\varepsilon)}^2\dd +\int_{\R^n}H_{ij}\partial_{ij}u_{(\xi,\varepsilon)} \bar w_{(\xi,\varepsilon)}\dd.
\end{split}
\end{equation}
In the estimates below, write $u=u_{(\xi,\varepsilon)}$ and $w=w_{(\xi,\varepsilon)}$, where $(\xi,\varepsilon)=\lambda(\bar\xi,\bar\varepsilon)$.

\begin{lemma}\label{lem:metric-energy}
There is a constant $C>0$, independent of $N$, such that, for all sufficiently large $N$,
\begin{equation}\label{eq:metric-energy-estimate}
\begin{aligned}
&\Bigl\|J_{\bar g}^{\R^n}(u)-J_\sph-\frac{1}{2}\int\bar h_{ik}\bar h_{jk}\partial_iu\partial_ju+\frac{c(n)}{4}\int(\partial_k\bar h_{ij})^2u^2\Bigr\|_{C^1(\overline\Omega)}\\
&\qquad\le C\Bigl[\lambda^{3(l+2)}(1+\log(\tfrac{\rho}{\lambda}))^{\mathbf1_{\{n=3l+8\}}}+(\lambda/\rho)^{n-2}\Bigr].
\end{aligned}
\end{equation}
\end{lemma}
\begin{proof}
The inverse-metric expansion on $B_\rho$ is
\[
\bar g^{ij}=\delta^{ij}-\bar h_{ij} +\frac{1}{2}\bar h_{ik}\bar h_{jk}+O(|\bar h|^3).
\]
Since $\bar h$ is divergence-free there, Lemma~\ref{lem:scalar-expansion} gives
\begin{equation*}
R_{\bar g}+\frac{1}{4}(\partial_k\bar h_{ij})^2 =O(|\bar h|^2|\partial^2\bar h|+|\bar h||\partial\bar h|^2).
\end{equation*}
The cubic term in the inverse metric and the scalar-curvature remainder multiplied by $u^2$ are both bounded on $B_\rho$ by the radial integral
\begin{equation}\label{eq:cubic-radial}
\begin{aligned}
&C\lambda^{n-2}\int_0^\rho(\lambda+r)^{3(l+2)+2-2n}r^{n-1}\,dr\\
&\qquad\le C\lambda^{3(l+2)}(1+\log(\frac{\rho}{\lambda}))^{\mathbf1_{\{n=3(l+2)+2\}}} +C\lambda^{n-2}.
\end{aligned}
\end{equation}
Inequality~\eqref{eq:cubic-radial} follows by splitting the integral at $r=\lambda$ and using $\lambda\le\lambda+r\le2\lambda$ for $r\le\lambda$ and $r\le\lambda+r\le2r$ for $r\ge\lambda$.

Contracting \eqref{eq:stress-identity} with the trace-free $\bar h$ and integrating twice by parts gives
\begin{equation}\label{eq:global-stress}
\int\bar h_{ij}\partial_iu\partial_ju =c(n)\int\partial_{ij}\bar h_{ij}u^2.
\end{equation}
Both integrations by parts are justified by the compact support of $\bar h$. Since $\bar h$ is divergence-free on $B_\rho$, the integrand on the right vanishes there. The contribution from $\rho<|x|<1$ is bounded in absolute value by $ C\lambda^{n-2}\int_\rho^1r^{3-n}dr \le C\rho^2(\lambda/\rho)^{n-2}. $

Since $\bar h$ is supported in $B_1$, the integrand in \eqref{eq:global-stress} also vanishes outside $B_1$. For the other exterior terms, bounded coefficients on $\rho<|x|<1$ give $ C\lambda^{n-2}\int_\rho^1(r^{1-n}+r^{3-n})\,dr\le C(\lambda/\rho)^{n-2}, $ and again the metric perturbation vanishes for $|x|>1$. Combining these estimates with \eqref{eq:cubic-radial} and \eqref{eq:global-stress} proves the $C^0$ part of \eqref{eq:metric-energy-estimate}.

Since $ |\partial_{\bar\theta}(\partial_iu\partial_ju)| \le C\lambda^{n-2}(\lambda+|x|)^{2-2n},\ \ |\partial_{\bar\theta}(u^2)|\le Cu^2, $ the same estimates give the $C^1$ bound.
\end{proof}
\begin{lemma}\label{lem:polynomial-replacement}
There is a constant $C>0$, independent of $N$, such that, for all sufficiently large $N$,
\begin{equation}\label{eq:polynomial-comparison}
\begin{aligned}
&\Bigl\|J_{\bar g}^{\R^n}(u)-J_\sph+\int H^\lambda_{ij}\partial_{ij}u\,w-\lambda^{2(l+2)}F(\bar\xi,\bar\varepsilon)\Bigr\|_{C^1(\overline\Omega)}\\
&\qquad\le C\Bigl[\lambda^{3(l+2)}(1+\log(\tfrac{\rho}{\lambda}))^{\mathbf1_{\{n=3l+8\}}}+(\lambda/\rho)^{n-2}\Bigr].
\end{aligned}
\end{equation}
The right-hand side of \eqref{eq:polynomial-comparison} is $o(\lambda^{2(l+2)})$ as $N\to\infty$.
\end{lemma}
\begin{proof}
By changing variables in \eqref{eq:auxiliary-function} and using \eqref{eq:exact-scaling}, the correction integral cancels from the difference estimated in \eqref{eq:polynomial-comparison}. Since $\bar h=H^\lambda$ on $B_\rho$, Lemma~\ref{lem:metric-energy} reduces the claim to
\begin{equation}\label{eq:model-replacement-exact}
\begin{aligned}
&\Biggl\|\frac{1}{2}\int_{\R^n\setminus B_\rho}(\bar h_{ik}\bar h_{jk}-H^\lambda_{ik}H^\lambda_{jk})\partial_iu\partial_ju\\
&\qquad-\frac{c(n)}{4}\int_{\R^n\setminus B_\rho}\bigl((\partial_k\bar h_{ij})^2-(\partial_kH^\lambda_{ij})^2\bigr)u^2\Biggr\|_{C^1(\overline\Omega)}\le C(\lambda/\rho)^{n-2}.
\end{aligned}
\end{equation}
On $B_1\setminus B_\rho$, the tensors $\bar h,H^\lambda$ and their first derivatives are uniformly bounded, and \eqref{eq:bubble-pointwise} gives $|\nabla u|^2+u^2\le C\lambda^{n-2}(\lambda+|x|)^{2(1-n)}$. On $\R^n\setminus B_1$, one has $\bar h=0$, $|H^\lambda(x)|\le C|x|^{l+2}$, $|\partial H^\lambda(x)|\le C|x|^{l+1}$, and $|\nabla u|^2+u^2\le C\lambda^{n-2}(\lambda+|x|)^{2(2-n)}$. Integration gives the $C^0$ bound in \eqref{eq:model-replacement-exact}; the exterior integral converges since $n>2l+8$.

The same bounds hold after applying $\partial_{\bar\theta}$ to the bubble factors, by \eqref{eq:bubble-pointwise}. This proves \eqref{eq:model-replacement-exact}, and Lemma~\ref{lem:metric-energy} gives \eqref{eq:polynomial-comparison}. Finally, division by $\lambda^{2(l+2)}$ gives a quantity tending to zero, since $n-2>2(l+2)$, $\lambda=2^{-N}$, and $\rho=(4N^2)^{-1}$.
\end{proof}

This identifies the limiting functional. Section~\ref{sec:global} will show that the actual reduced functional, evaluated on the projected family, has the same leading term.

\section{Endpoint critical points}\label{sec:auxiliary}
We choose the coefficients $\tau_{a,b}$ and $\beta$ in $H$ so that $F$ has a nondegenerate critical point at $(0,1)$ with negative value. Proposition~\ref{prop:auxiliary-symmetry} gives translation stationarity, Proposition~\ref{prop:criticality-polynomials} imposes scale criticality, and Proposition~\ref{prop:endpoint-critical} establishes the required Hessian signature. Throughout this section we drop the bars from the normalized parameters.

Put $\tau:=\tau_{0,0}$. For $l=2$, $n\ge52$, use
\begin{equation}\label{eq:model2}
H_{ij}=(\tau+\varrho^2)\Phi(W)_{ij}+\beta t\Phi(V)_{ij}.
\end{equation}
For $l=4$, $30\le n\le51$, use
\begin{equation}\label{eq:model4}
H_{ij}=(\tau-7\varrho^2+\varrho^4)\Phi(W)_{ij}+\beta t\Phi(V)_{ij}.
\end{equation}
For $l=6$, $25\le n\le29$, use
\begin{equation}\label{eq:model6}
H_{ij}=(\tau+48t^2+95\varrho^2-5t^4-12t^2\varrho^2 -19\varrho^4+\varrho^6)\Phi(W)_{ij} +\beta t\Phi(V)_{ij}.
\end{equation}

For $l=2,4$, the odd $V$-term removes an axial translation degeneracy. For $l=6$, the even $t$-dependent terms suffice, and the successful choice in Proposition~\ref{prop:endpoint-critical} has $\beta=0$. 

\subsection{Evaluation and coefficients}
For $i>-1$ and $2j-i>1$, define $\mci_j^i:=\int_0^\infty r^i(1+r^2)^{-j}dr, $ and set  $\mu_n:=\frac{3|S^{n-1}|}{2n(n+2)}$.

The following identities follow from the calculations in \cite[proof of Proposition~19]{Brendle} and \cite[proof of Proposition~21]{Brendle}, respectively,  using $H_{ii}=0$ and $x_iH_{ij}=0$.

\begin{lemma}[Brendle \cite{Brendle}]
For $F$ defined by \eqref{eq:auxiliary-function}, with indices $1\le p,q\le n$,
\begin{equation}\label{eq:F-axis-value}
F(0,\varepsilon)=-\frac{c(n)}{4}\varepsilon^{n-2} \int_{\R^n}(\varepsilon^2+|x|^2)^{2-n}|\partial H|^2\dd,
\end{equation}
and
\begin{equation}\label{eq:F-hessian}
\begin{split}
\partial_{\xi_p\xi_q}^2F(0,\varepsilon)={}& (n-2)^2\varepsilon^{n-2}\int (\varepsilon^2+|x|^2)^{-n}H_{pk}H_{qk}\\
&-\frac{(n-2)^2}{4}\varepsilon^{n-2}\int (\varepsilon^2+|x|^2)^{-n}|\partial H|^2x_px_q\\
&+\frac{(n-2)c(n)}{2}\varepsilon^{n-2}\int (\varepsilon^2+|x|^2)^{1-n}|\partial H|^2\delta_{pq}.
\end{split}
\end{equation}
\end{lemma}
\begin{proposition}[{Translation stationarity and Hessian block structure}]\label{prop:auxiliary-symmetry}
For all parameters,
\begin{equation}\label{eq:reflection}
F(\xi_1,\xi',\varepsilon)=F(\xi_1,-\xi',\varepsilon).
\end{equation}
Hence $\partial_{\xi'}F(0, \varepsilon)=0$ and   all mixed Hessian blocks between $\xi_1$, $\xi'$, and $\varepsilon$ vanish at $(0,1)$. Moreover, $\partial_{\xi_1}F(0,\varepsilon)=0$ for all $\varepsilon>0$.
\end{proposition}
\begin{proof}
The tensor $H$ is even under $x'\mapsto-x'$ and satisfies $H_{1j}=0$. The first two terms of \eqref{eq:auxiliary-function} are invariant under this reflection, and the linear correction equation and constraints are preserved. Uniqueness gives the same invariance for $\bar w$ and its contribution to $F$, proving \eqref{eq:reflection}. It follows that $\partial_{\xi_j}F$, $\partial_{\xi_1\xi_j}^2F$, and $\partial_{\varepsilon\xi_j}^2F$ vanish when $\xi'=0$, for $j\ge2$.

Now we prove $\partial_{\xi_1}F(0,\varepsilon)=0$. On the distinguished axis $\xi=se_1$, the identities $H(x)x=0$ and $H_{1j}=0$ imply $H_{ij}\partial_{ij}u_{(se_1,\varepsilon)}=0$ and $H_{ik}H_{jk}\partial_iu\partial_ju=0$. Thus $\bar w_{(se_1,\varepsilon)}=0$ and
\[
F(se_1,\varepsilon)=-\frac{c(n)}{4}\int|\partial H|^2u_{(se_1,\varepsilon)}^2\dd.
\]
Differentiation at $s=0$ gives
\begin{equation}\label{eq:xi1-derivative}
\partial_{\xi_1}F(0,\varepsilon) =-\frac{(n-2)c(n)}{2}\varepsilon^{n-2} \int t|\partial H|^2(\varepsilon^2+|x|^2)^{1-n}\dd.
\end{equation}
Write $H=H_W+H_V$ with $H_V:=\beta t\Phi(V)$. The terms $|\partial H_W|^2$ and $|\partial H_V|^2$ are even in $t$, so they give zero in \eqref{eq:xi1-derivative}. For the cross term, at each fixed $t,\varrho$,  the primed angular integrals reduce to
\[
\int_{S^{n-2}}\ip{\Phi(W)}{\Phi(V)}=0,\qquad \int_{S^{n-2}}\sum_{\alpha=2}^n \ip{\partial_\alpha\Phi(W)}{\partial_\alpha\Phi(V)}=0.
\]
These follow by polarization of the sphere identities in Lemma~\ref{lem:weyl-contractions} and $\ip WV=0$. Terms containing $x'\cdot\nabla\Phi(V)$ reduce to twice the first integral by homogeneity. Hence the cross term also vanishes, proving $\partial_{\xi_1}F(0,\varepsilon)=0$. Differentiating this identity in $\varepsilon$ gives the last mixed-block vanishing.
\end{proof}

For \eqref{eq:model2} and \eqref{eq:model4} with $\beta=0$, the tensor $H$ is independent of $t$. The formula for $F(se_1,\varepsilon)$ in the proof of Proposition~\ref{prop:auxiliary-symmetry} makes $F(se_1,\varepsilon)$ independent of $s$, by a change of variable in the $t$ integral. Thus the second derivative in the $\xi_1$ direction vanishes. For $H_V:=\beta t\Phi(V)$, the pure $V$ contribution contains a negative multiple of
\[
\int\bigl(|\Phi(V)|^2+t^2|\nabla\Phi(V)|^2\bigr) u_{(se_1,\varepsilon)}^2.
\]
After translating $t$, its second derivative in $s$ is strictly negative for $\beta\ne0$. The angular orthogonality implies that the first derivative and mixed Hessian entries vanish. For \eqref{eq:model6}, the even $t$-dependent terms already remove the degeneracy, so the odd component will be set to zero.

Set $\tau_{0,l/2}=1$ and the other coefficients with $a+b=l/2$ equal to zero. Then
\[
H=\sum_{a+b\le l/2}\tau_{a,b}H^{a,b}+\beta t\Phi(V).
\]
Recall that $H^{a,b}=t^{2a}\varrho^{2b}\Phi(W)$.

The computation proceeds from the following angular integrals $I,A_{pq},B_{pq}$ to the sums $Q_l,M^l$, and then to their coefficients $\qhat_s^l,X_l,Y_l,Z_l$. Propositions~\ref{prop:criticality-polynomials} and \ref{prop:endpoint-critical} use these coefficients to impose criticality and establish nondegeneracy.

We first record the angular contractions of these homogeneous terms. Define $$I(a,b,c,d):=\int_{S^{n-1}}\partial_kH^{a,b}_{ij}\partial_kH^{c,d}_{ij},
$$
$$A_{pq}(a,b,c,d):=\int_{S^{n-1}}H^{a,b}_{pi}H^{c,d}_{iq},\ \ \ B_{pq}(a,b,c,d):=\int_{S^{n-1}}\partial_kH^{a,b}_{ij}\partial_kH^{c,d}_{ij}x_px_q.
$$
Define $Q_l$ and $M^l=(M^l_{pq})$ by
\begin{align}
Q_l(\tau,\varepsilon)&:=\sum\tau_{a,b}\tau_{c,d}I(a,b,c,d)\mci_{n-2}^{2(a+b+c+d)+n+1}\varepsilon^{2(a+b+c+d)+4},\label{eq:Q-sum}\\
M^l_{pq}(\tau,\varepsilon)&:=\sum\tau_{a,b}\tau_{c,d} \Bigl(4A_{pq}-B_{pq}+\frac{I\delta_{pq}}{n+2(a+b+c+d)+2}\Bigr)\nonumber\\
&\qquad\times\mci_n^{2(a+b+c+d)+n+3}\varepsilon^{2(a+b+c+d)+2}.\label{eq:M-sum}
\end{align}
Both sums are over nonnegative integers $a,b,c,d$ with $a+b\le l/2$ and $c+d\le l/2$. The arguments $(a,b,c,d)$ of $I,A_{pq},B_{pq}$  are suppressed in \eqref{eq:M-sum}. Proposition~\ref{prop:angular-contractions} expresses these integrals in terms of the beta-function ratios \eqref{eq:Rratio}. We now collect the scalar polynomial coefficients and the three coefficients of the translation matrix. Define $\qhat_s^l$, $X_l,Y_l,Z_l$ by
\begin{align}
Q_l(\tau,\varepsilon)&=\mu_n|W|^2\mci_{n-2}^{n+1}\varepsilon^4 \sum_{s=0}^l\qhat_s^l(\tau)\varepsilon^{2s},\nonumber\\
M^l_{pq}(\tau,1)&=\mu_n\mci_n^{n+5}\Bigl[X_l(\tau)(W\circ W)_{pq}\label{eq:XYZ-definition}\\
&\quad+|W|^2\bigl(Y_l(\tau)(\delta_{pq}-\delta_{p1}\delta_{q1}) +Z_l(\tau)\delta_{p1}\delta_{q1}\bigr)\Bigr].\nonumber
\end{align}
This normalization gives $\int_{S^{n-1}}|\Phi(W)|^2=\mu_n|W|^2$. Their formulas are listed in Appendix~\ref{sec:coefficients}.

The coefficients now determine both the value of $F$ on the scale axis and its Hessian. We first evaluate $F(0,\varepsilon)$ and determine the parameters for which $(0,1)$ is a critical point.

Define the polynomials
\begin{equation}
P_2(\tau):=(n-10)(n-8)\tau^2+3(n-10)(n+3)\tau +2(n+3)(n+7),\label{eq:P2}
\end{equation}
\begin{align}
P_4(\tau):={}&(n-14)(n-12)(n-10)(n-8)\tau^2\nonumber\\
&-(n+3)(n-14)(n-12)(17n-230)\tau\nonumber\\
&+(n+3)(66n^3-1845n^2+5260n+143983),\label{eq:P4}\\
P_6(\tau):={}&(n-18)(n-16)(n-14)(n-12)(n-10)(n-8)\tau^2\nonumber\\
&+(n-18)(n-16)(n-14)\nonumber\\
&\quad\times(214n^3-4940n^2+19105n+136533)\tau\nonumber\\
&+10549n^6-639262n^5+12735837n^4-53498987n^3\nonumber\\
&-1021178224n^2+7036933380n+30233004948.\label{eq:P6}
\end{align}
By Lemma~\ref{lem:endpoint-coefficient-signs}, each $P_l$ has two distinct real roots, $\tau_{l,-}<\tau_{l,+}$, in the following respective ranges: $l=2$ and $n\ge 52$; $l=4$ and $30\le n\le 51$;  $l=6$ and $25\le n\le 29$.
\begin{proposition}[Evaluation and criticality]\label{prop:criticality-polynomials}
For $H$ given by \eqref{eq:model2}--\eqref{eq:model6},
\begin{equation}\label{eq:F-coefficient-value}
F(0,\varepsilon)=-\frac{c(n)}{4}\mu_n|W|^2\mci_{n-2}^{n+1} \Bigl(\varepsilon^4\sum_{s=0}^l\qhat_s^l(\tau)\varepsilon^{2s} +3\beta^2\frac{n+2}{n-8}\varepsilon^6\Bigr).
\end{equation}
Then $(0,1)$ is a critical point of $F$ if and only if
\begin{equation}\label{BBB}
\beta^2= \left\{
\begin{array}{lr}
-\frac{4P_2(\tau)}{9(n-10)},& \text{if}\ l=2,\\
-\frac{4P_4(\tau)}{9(n-14)(n-12)(n-10)}, &  \text{if}\ l=4,\\
-\frac{4P_6(\tau)}{9(n-18)(n-16)(n-14)(n-12)(n-10)}, & \text{if}\ l =6.
\end{array}
\right.
\end{equation}\end{proposition}
\begin{proof}
Write $H=H_W+\beta t\Phi(V)$. The mixed term in $|\partial H|^2$ has zero angular integral by the orthogonality used in Proposition~\ref{prop:auxiliary-symmetry}. Separating radial and angular variables in \eqref{eq:F-axis-value} therefore gives the $W$ contribution $-c(n)Q_l/4$. For the odd term, Proposition~\ref{prop:odd-tensor-contribution} gives $I_V=3\mu_n\beta^2|V|^2$ and a radial factor $\mci_{n-2}^{n+3}\varepsilon^6$. Using $|V|=|W|$ and $\mci_{n-2}^{n+3}/\mci_{n-2}^{n+1}=(n+2)/(n-8)$ proves \eqref{eq:F-coefficient-value}.

We spell out the coefficient calculation for $l=2$. Here $H_W=(\tau+\varrho^2)\Phi(W)$, and Proposition~\ref{prop:angular-contractions} gives $\frac{I(0,0,0,0)}{\mu_n|W|^2}=2(n+2)$, $\frac{I(0,0,0,1)}{\mu_n|W|^2}=2(n+3)$, and $\frac{I(0,1,0,1)}{\mu_n|W|^2}=2(n+3)(n+7)/(n+4)$. For the middle identity, use $2(n+2)R_1(0,1)+4R_2(0,0)$ with $R_1(0,1)=(n+1)/(n+2)$; in the last identity, $bd+b+d=3$. The two mixed terms in \eqref{eq:Q-sum} have the same coefficient $\tau$. Together with the radial ratios $\mci_{n-2}^{n+3}/\mci_{n-2}^{n+1}=(n+2)/(n-8)$ and $\mci_{n-2}^{n+5}/\mci_{n-2}^{n+1}=(n+2)(n+4)/((n-8)(n-10))$, this yields $\qhat_0^2=2(n+2)\tau^2$, $\qhat_1^2=4(n+2)(n+3)\tau/(n-8)$, and $\qhat_2^2=2(n+2)(n+3)(n+7)/((n-8)(n-10))$. The same bilinear sums give the lists for $l=4,6$ in Appendix~\ref{sec:coefficients}, with the powers of $t$ retained by $R_r(a+c,b+d)$.

Differentiating \eqref{eq:F-coefficient-value} at $\varepsilon=1$ shows that the scale derivative vanishes precisely when
\[
\sum_{s=0}^l(2s+4)\qhat_s^l(\tau) +18\beta^2\frac{n+2}{n-8}=0.
\]
Substitution of the three scalar coefficient formulas and collection of powers of $\tau$ give
\[
\sum_{s=0}^l(2s+4)\qhat_s^l(\tau) =\frac{8(n+2)P_l(\tau)}{\prod_{j=0}^{l-1}(n-8-2j)}.
\]
Solving for $\beta^2$ proves \eqref{BBB}.
\end{proof}

\subsection{The full Hessian}
The preceding proposition imposes criticality. The next proposition computes the Hessian at that critical point and chooses the coefficients so that it is nondegenerate.

\begin{proposition}[{Nondegenerate  critical points}]\label{prop:endpoint-critical}
For each endpoint pair $(n,l)$, the coefficients in \eqref{eq:model2}--\eqref{eq:model6} can be chosen so that
\[
\nabla F(0,1)=0,\qquad F(0,1)<0,\qquad \det D^2F(0,1)\ne0.
\]
For $l=2,4$, one can take $\beta>0$, and the Hessian has exactly one negative eigenvalue. For $l=6$, the choice $\tau=\tau_{6,-}$ and $\beta=0$ gives a positive definite Hessian and hence a strict local minimum. 

More generally, whenever $\tau$ and $\beta\ge0$ satisfy \eqref{BBB}, the Hessian is given by the following formulas. At $(0,1)$, all mixed Hessian blocks between $\xi_1$, $\xi'$, and $\varepsilon$ vanish. For $2\le p,q\le n$, the transverse block is
\begin{equation}\label{eq:transverse-Hessian}
\begin{split}
\partial_{\xi_p\xi_q}^2F(0,1) =&\frac{(n-2)^2\mu_n\mci_n^{n+5}}{4}\\
&\Bigl[\bigl(X_l(\tau)-\frac{4\beta^2}{n+4}\bigr)(W\circ W)_{pq}+|W|^2Y_l(\tau)\delta_{pq}\Bigr].
\end{split}
\end{equation}
The $\xi_1$ and $\varepsilon$ diagonal entries are
\begin{equation}\label{eq:distinguished-Hessian}
\partial_{\xi_1\xi_1}^2F(0,1) =\frac{(n-2)^2\mu_n\mci_n^{n+5}}{4}|W|^2 \Bigl(Z_l(\tau)-\frac{4\beta^2}{n+4}\Bigr).
\end{equation}
\begin{equation}\label{eq:scale-Hessian}
\partial_{\varepsilon\varepsilon}^2F(0,1) =-\frac{c(n)}{4}\mu_n|W|^2\mci_{n-2}^{n+1} \sum_{s=0}^l(2s+4)(2s-2)\qhat_s^l(\tau).
\end{equation}
\end{proposition}
\begin{proof}
We first compute the Hessian, and then choose the parameters. For the $W$ terms, \eqref{eq:F-hessian} and \eqref{eq:M-sum} give the translation matrix $(n-2)^2M^l(\tau,1)/4$. Substituting \eqref{eq:XYZ-definition} identifies its transverse and distinguished blocks. For the odd term, use the contractions in Proposition~\ref{prop:odd-tensor-contribution}. The $W$--$V$ cross terms in the transverse and distinguished Hessian blocks vanish by parity in $t$. Use the radial identity $\mci_{n-1}^{n+3}=\frac{2(n-1)}{n+4}\mci_n^{n+5}$. For $p,q>1$, the transverse contribution is determined by
\begin{equation*}
4(A_V)_{pq}-(B_V)_{pq}+\frac{I_V\delta_{pq}}{n+4}=-\frac{4\mu_n\beta^2}{n+4}(V\circ V)_{pq}.
\end{equation*}
Thus the scalar $\delta_{pq}$ terms cancel, and the coefficient of $(V\circ V)_{pq}$ comes from $4-8=-4$. In the distinguished direction, $(A_V)_{11}=0$, so
\[
4(A_V)_{11}-(B_V)_{11}+\frac{I_V}{n+4}=-\frac{4\mu_n\beta^2|V|^2}{n+4}.
\]
The coefficient of $|V|^2$ comes from $-7+3=-4$. Multiplying both expressions by the common radial factor $(n-2)^2\mci_n^{n+5}/4$ in \eqref{eq:F-hessian}, and using $V\circ V=W\circ W$ and $|V|=|W|$, gives the $V$ contributions in \eqref{eq:transverse-Hessian} and \eqref{eq:distinguished-Hessian}.

For the scale entry, differentiate \eqref{eq:F-coefficient-value} twice. At a critical point, the first derivative gives $18\beta^2(n+2)/(n-8)=-\sum_s(2s+4)\qhat_s^l$. Consequently the coefficient in the second derivative is
\[
\sum_s(2s+4)(2s+3)\qhat_s^l -5\sum_s(2s+4)\qhat_s^l =\sum_s(2s+4)(2s-2)\qhat_s^l,
\]
which proves \eqref{eq:scale-Hessian}. The mixed blocks vanish by Proposition~\ref{prop:auxiliary-symmetry}.

It remains to make these blocks nondegenerate. The matrix $W\circ W$ is nonnegative and has trace $|W|^2$, so $0\le W\circ W\le |W|^2\Id$. Hence
\[
X_l(W\circ W)+Y_l|W|^2\Id \ge\bigl(Y_l+\min(X_l,0)\bigr)|W|^2\Id.
\]
Lemma~\ref{lem:endpoint-coefficient-signs} makes the right-hand coefficient positive at the specified roots and makes the sum in \eqref{eq:scale-Hessian} negative there. Thus the transverse block and scale entry are positive when $\beta=0$. For $l=6$, the same lemma gives $Z_6(\tau_{6,-})>0$, so the full Hessian is positive definite.

For $l=2$, choose $\tau$ just to the right of $\tau_{2,-}$; for $l=4$, choose it just to the left of $\tau_{4,+}$. Then $P_l(\tau)<0$, and the corresponding case of \eqref{BBB} determines $\beta>0$, tending to zero as $\tau$ approaches the root. Continuity preserves positivity of the transverse and scale blocks. Since $Z_2=Z_4=0$, \eqref{eq:distinguished-Hessian} is strictly negative. The Hessian therefore has exactly one negative direction.

Fix these choices of $\tau,\beta$ for all $N$. A common rescaling of $W,V$ multiplies $F$ and its Hessian by the square of the rescaling factor, preserving criticality and signature. After this rescaling is fixed, the absolute values of the Hessian eigenvalues have a positive lower bound independent of $N$. Finally, $F(0,1)<0$ follows from \eqref{eq:F-axis-value}, because $H$ contains the nonzero highest-degree term $T_\glo$ and hence $\partial H$ is not identically zero.
\end{proof}

\section{Global construction and finite-dimensional reduction}\label{sec:global}
Recall from Section~\ref{sec:euclidean-estimates} that $m=(n-2)/2$, $p=(n+2)/(n-2)$, and $q=2n/(n+2)$.
We first extend the local construction to a single metric on $M$ with positive Yamabe constant and the required global Weyl nonvanishing. For this fixed metric, a uniform projected inverse yields a projected family whose normalized reduced functional converges to $F$ in $C^1$. The critical point of $F$ persists and yields solutions. We first treat the endpoint orders and then the higher prescribed orders.

\subsection{Geometric preparations}
We first record the known vanishing obstruction that explains the sharpness of Theorem~\ref{thm:global-blowup}; the metric construction begins with Lemma~\ref{lem:positive-flat}. The following vanishing statements are consequences of the Khuri--Marques--Schoen estimates, as recorded in \cite[p.~3 of the preprint]{MarquesBlowup}.

\begin{proposition}[Khuri--Marques--Schoen; Marques]\label{prop:fixed-order-vanishing}
Let $p$ be a blow-up point for a sequence of positive solutions of the Yamabe equation on a fixed smooth closed manifold. Then
\begin{align*}
\nabla_g^2W_g(p)=\nabla_g^3W_g(p)&=0&&\text{if }25\le n\le51,\\
\nabla_g^4W_g(p)=\nabla_g^5W_g(p)&=0&&\text{if }25\le n\le29.
\end{align*}
The lower derivatives vanish as in the classical results $W_g(p)=0$ for $n\ge6$ and $\nabla_gW_g(p)=0$ for $n\ge8$.
\end{proposition}
On a closed manifold, positivity of $Y(M,[g])$ implies coercivity of $Q_g$ in $H^1(M)$. Indeed the definition of $Y$ gives a positive lower bound in the critical $L^{2n/(n-2)}$ norm, hence in $L^2$ on a compact manifold. Since $R_g$ is bounded, the identity for $Q_g$ then also controls $\|\nabla f\|_2^2$.

\begin{lemma}[{A positive Yamabe type manifold  with a Euclidean region}]\label{lem:positive-flat}
Let $M^n$, $n\ge3$, be closed and of positive Yamabe type, and fix $p\in M$. There is a smooth metric $g_*$ on $M$ for which $Q_{g_*}$ is coercive and which is exactly Euclidean on a coordinate ball $B_3$ centered at $p$.
\end{lemma}
\begin{proof}
Start with any smooth metric $g'$ for which $Y(M,[g'])>0$, and use normal coordinates $y$ at $p$. For small $r>0$, choose a radial cutoff $\zeta_r$, equal to one on $B_r$ and zero outside $B_{2r}$, with $|\partial^k\zeta_r|\le C_kr^{-k}$ for $k\le2$, and set $ g_r:=g'+\zeta_r(\delta-g') $ inside the coordinate patch and $g_r=g'$ outside. In these coordinates $ g'-\delta=O(|y|^2), \partial g'=O(|y|), \partial^2g'=O(1). $ Consequently $g_r$ is positive definite for small $r$, the metric and volume-density differences are $O(r^2)$, and the first two coordinate derivatives of $g_r$ are uniformly bounded. In particular $R_{g_r}$ is uniformly bounded on $B_{2r}$.

Using the volume form of $g'$ to compare the two quadratic forms gives
\begin{equation*}
|Q_{g_r}(f)-Q_{g'}(f)| \le Cr^2\|\nabla_{g'}f\|_2^2+ C\|\mathbf1_{B_{2r}}\|_{L^{n/2}}\|f\|_{L^{\frac{2n}{n-2}}}^2 \le Cr^2\|f\|_{H^1(g')}^2.
\end{equation*}
The factor $r^2$ in the potential term is the $L^{n/2}$ norm of a bounded function supported on a ball of volume $O(r^n)$. Sobolev and coercivity of $Q_{g'}$ imply positivity and coercivity of $Q_{g_r}$ for sufficiently small $r$.

Fix $r>0$ sufficiently small that $Q_{g_r}$ is coercive, and set $a:=r/3$. In the rescaled coordinates $x:=a^{-1}y$, the ball $B_r$ becomes $B_3$, and $ g_r=a^2\delta_{ij}\,dx^i\,dx^j\ \text{on }B_3. $ Define $g_*:=a^{-2}g_r$. Then $g_*=\delta$ on $B_3$. For every constant $c>0$, $ Q_{cg}(f)=c^{\frac{n-2}{2}}Q_g(f). $ The $H^1$ norms of two metrics differing by a fixed positive factor are equivalent, so this rescaling preserves coercivity. Thus $g_*$ has the required properties.
\end{proof}

\begin{proposition}[{Global Weyl-derivative nonvanishing by small perturbations}]\label{prop:weyl-generic}
Let $(M^n,g)$ be a closed smooth Riemannian manifold with $n\ge4$, let $l$ be a positive integer, and let $h$ be a smooth reference metric on $M$. For every $\epsilon>0$, there exists a smooth metric $\widetilde g$ such that
\[
\|\widetilde g-g\|_{C^{l+2}(h)}<\epsilon, \qquad |\nabla_{\widetilde g}^lW_{\widetilde g}|_{\widetilde g}>0 \quad\text{on }M.
\]
\end{proposition}
\begin{proof}
Let
\[
S_0^2:=\{A\in\operatorname{Sym}_{n-1}(\mathbb R):\operatorname{tr}A=0\}
\]
be the space of trace-free symmetric $(n-1)\times(n-1)$ matrices. Its dimension is $ r:= -1+n(n-1)/2>n. $ For a coordinate chart $U$, let $J^{l+2}\operatorname{Met}(U)$ denote the bundle of metric jets of order $l+2$. A jet $j_x^{l+2}q$ records the metric coefficients and their coordinate derivatives through order $l+2$ at $x$. Applying Gram--Schmidt to the coordinate frame gives an orthonormal frame $e_1,\ldots,e_n$ depending smoothly on  $q$. Define $ F_U:J^{l+2}\operatorname{Met}(U)\longrightarrow S_0^2 $ by
\[
(F_U(j_x^{l+2}q))_{ab} := (\nabla_q^lW_q)_x(e_1^{\otimes l};e_1,e_a,e_1,e_b),
\]
where $2\le a,b\le n$ and the first $l$ slots are derivative slots. Covariant differentiation preserves the Weyl symmetries and trace identities, so this matrix is symmetric and trace-free.

Fix a metric jet at $p$ and a matrix $A\in S_0^2$.

Choose $q$-normal coordinates at $p$ with $\partial_i|_p=e_i|_p$, and set $v_{ab}:=(x^1)^{l+2}A_{ab}/(l+2)!$ for $2\le a,b\le n$ and $v_{1i}=0$. Since $j_p^{l+1}v=0$ and $\partial_1^{l+2}v_{ab}(p)=A_{ab}$, only the principal curvature terms contribute to $DF_U$. The variations of $\partial_1^l\Rm_{1a1b}$ and $\partial_1^l\Ric_{ab}$ are both $-A_{ab}/2$; those of $\partial_1^l\Ric_{11}$ and $\partial_1^lR$ vanish because $\operatorname{tr}A=0$. Thus the Weyl decomposition gives

\begin{equation}\label{eq:weyl-jet-submersion}
DF_U|_{j_p^{l+2}q}[j_p^{l+2}v] =-\frac{n-3}{2(n-2)}A.
\end{equation}
Hence the restriction of $F_U$ to each metric-jet fiber is a submersion.

Choose finitely many coordinate neighborhoods $U_i\Subset V_i$ such that the $U_i$ cover $M$, and choose $\chi_i\in C_c^\infty(V_i)$ equal to one near $\overline{U_i}$. In the coordinates on $V_i$, multiply $\chi_i$ by all symmetric covariant two-tensors whose coefficients are monomials of degree at most $l+2$, and extend the resulting tensors by zero. Enumerate them as $s_1,\ldots,s_J$. At every $p\in U_i$, the tensors from this chart realize arbitrary Taylor data through order $l+2$. Consequently
\[
\operatorname{span}\{j_p^{l+2}s_1,\ldots,j_p^{l+2}s_J\} =J_p^{l+2}(S^2T^*M) \qquad(p\in M).
\]
For $t\in B_\delta(0)\subset\mathbb R^J$, set
\[
g_t:=g+\sum_{\alpha=1}^Jt_\alpha s_\alpha.
\]
Choose $\delta>0$ so that every $g_t$ is positive definite. For each $i$, define
\[
\Phi_i:B_\delta(0)\times U_i\longrightarrow S_0^2, \qquad \Phi_i(t,p):=F_{U_i}(j_p^{l+2}g_t).
\]
At every $(t,p)$, the parameter derivative of $j_p^{l+2}g_t$ spans the full tangent space of the metric-jet fiber. Together with \eqref{eq:weyl-jet-submersion}, this proves that $D_t\Phi_i$ is onto. This uses the same parameter family on all charts, including their overlaps.

By the regular value theorem, $Z_i:=\Phi_i^{-1}(0)$ is a smooth submanifold of dimension $J+n-r<J$. The coordinate projection $(t,p)\mapsto t$ from $Z_i$ to $B_\delta(0)$ has rank less than $J$ everywhere. Sard's theorem implies that its image has measure zero; see \cite[Chapter~II]{GolubitskyGuillemin}. The finite union of these images also has measure zero. Choose $t$ outside this union arbitrarily close to zero. Then $\Phi_i(t,p)\ne0$ for every $p\in U_i$ and every $i$, which implies $\nabla_{g_t}^lW_{g_t}\ne0$ on $M$. Finally,
\[
\|g_t-g\|_{C^{l+2}(h)} \le\sum_{\alpha=1}^J|t_\alpha|\,\|s_\alpha\|_{C^{l+2}(h)} \longrightarrow0 \qquad\text{as }t\to0.
\]
Taking $t$ sufficiently small and setting $\widetilde g:=g_t$ proves the proposition.
\end{proof}

\begin{lemma}[{Global completion preserving the local metric}]\label{lem:metric-extension}
Let $g_*=\delta$ on $B_3$ be as in Lemma~\ref{lem:positive-flat}. Suppose that $h_0$ is smooth, symmetric, trace-free, supported in $B_1$, sufficiently small in $C^2$, and that $g_0=\exp(h_0)$ satisfies
\[
\nabla_{g_0}^jW_{g_0}(0)=0\quad(j<l),\qquad \nabla_{g_0}^lW_{g_0}(0)\ne0.
\]
There exist a smooth metric with positive Yamabe constant $\widehat g$ on $M$ and a fixed $b\in (0, 1/4)$ such that $ \widehat g=g_0\text{ on }B_{3b} $ and
\begin{equation}\label{eq:global-geometric-condition}
|\nabla_{\widehat g}^lW_{\widehat g}|>0 \quad\text{on }M.
\end{equation}
\end{lemma}
\begin{proof}
Insert $g_0$ in $B_1$ and leave $g_*$ unchanged outside $B_1$. The resulting metric $g_1$ is smooth, since $h_0$ vanishes near $\partial B_1$. A sufficiently small $C^2$ perturbation of a fixed smooth metric changes $Q_g(f)$ by at most $C\|h_0\|_{C^2}\|f\|_{H^1}^2$. Thus $Q_{g_1}$ remains coercive.

Since $\nabla_{g_1}^lW_{g_1}(0)\ne0$, choose $b\in (0, 1/4)$  and $\eta>0$ such that $ |\nabla_{g_1}^lW_{g_1}|_{g_1}\ge2\eta \ \text{on }\overline B_{4b}. $ Choose a smooth function $0\le\chi\le1$ equal to zero on $B_{3b}$ and to one on $M\setminus B_{7b/2}$. By Proposition~\ref{prop:weyl-generic}, for every $\epsilon>0$ there is a smooth metric $\widetilde g$ with $ \|\widetilde g-g_1\|_{C^{l+2}(g_1)}<\epsilon, \  |\nabla_{\widetilde g}^lW_{\widetilde g}|_{\widetilde g}>0 \ \text{on }M. $

Set $\widehat g:=g_1+\chi(\widetilde g-g_1)$. This is a Riemannian metric satisfying
\[
\widehat g=g_0\quad\text{on }B_{3b}, \qquad \widehat g=\widetilde g\quad\text{on }M\setminus B_{7b/2},
\]
and, for a constant depending only on the fixed cutoff, $g_1$, and $l$,
\[
\|\widehat g-g_1\|_{C^{l+2}(g_1)}\le C\|\widetilde g-g_1\|_{C^{l+2}(g_1)}<C\epsilon.
\]
For sufficiently small $\epsilon$, continuity preserves the lower bound $ |\nabla_{\widehat g}^lW_{\widehat g}|_{\widehat g}\ge\eta \ \text{on }\overline B_{4b}, $ and coercivity of $Q_{\widehat g}$. On $M\setminus B_{7b/2}$ the Weyl derivative is nonzero because $\widehat g=\widetilde g$. These sets cover $M$, proving \eqref{eq:global-geometric-condition}. Since $g_0=\exp(h_0)$ and $\tr h_0=0$, one also has $dV_{\widehat g}=dx$ on $B_{3b}$.
\end{proof}

At an endpoint, Lemma~\ref{lem:smooth-bump-series} and Lemma~\ref{lem:weyl-vanishing-order} verify the hypotheses with the original $h_0$ of \eqref{eq:h0}. For higher $l$, the analogous verification is given in Proposition~\ref{prop:higher-vanishing-orders}. Choose the tensor amplitudes small enough to apply Lemma~\ref{lem:metric-extension}, and fix the resulting metric $\widehat g$ and radius $b$ for all $N$. Both constructions give
\begin{equation}\label{eq:compact-coercivity}
c_0\|f\|_{H^1(M)}^2\le Q_{\widehat g}(f) \le C_0\|f\|_{H^1(M)}^2, \qquad \|f\|_{L^{\frac{2n}{n-2}}}\le C_0Q_{\widehat g}(f)^{1/2}.
\end{equation}
We use the $H^1(M,\widehat g)$ norm, and $H^{-1}(M)$ denotes its dual.

\subsection{The projected inverse}\label{subsec:projected-inverse}
We first treat the endpoint pairs $l=2$, $n\ge52$; $l=4$, $30\le n\le51$; and $l=6$, $25\le n\le29$. Apply Lemma~\ref{lem:metric-extension} to $g_0=\exp(h_0)$, where $h_0$ is the tensor constructed in Section~2, and fix the resulting metric $\widehat g$.

Fix $0\le\chi\le1$, smooth, equal to one on $B_b$ and supported in $B_{2b}$. We regard chart functions multiplied by $\chi$ as functions on $M$, extended by zero. For large $N$, the centers $x_N$ and the balls $B_{\rho_N}(x_N)$ lie in $B_{b/4}$. Continue to use
\[
(\xi,\varepsilon)=\lambda(\bar\xi,\bar\varepsilon),\qquad \lambda=\lambda_N,\qquad\rho=\rho_N.
\]
On the chart we evaluate the original $u_{(\xi,\varepsilon)}$, $w_{(\xi,\varepsilon)}$, $\varphi_{(\xi,\varepsilon,j)}$, and $Z_{(\xi,\varepsilon,j)}$ at $y=x-x_N$.

The normalized parameter derivatives of $\chi u_{(\xi,\varepsilon)}$ are
\begin{equation*}
\partial_{\bar\varepsilon}(\chi u_{(\xi,\varepsilon)})=\chi\lambda Z_{(\xi,\varepsilon,0)},\qquad \partial_{\bar\xi_j}(\chi u_{(\xi,\varepsilon)})=\chi\lambda Z_{(\xi,\varepsilon,j)},\quad 1\le j\le n.
\end{equation*}
Define
\begin{equation}\label{eq:compact-complement}
\E_{(\xi,\varepsilon)}(M):= \left\{f\in H^1(M):\int_M\chi\varphi_{(\xi,\varepsilon,j)} f\,dV_{\widehat g}=0, \quad0\le j\le n\right\}.
\end{equation}
In proofs below, we write $\E(M)$ when $(\xi,\varepsilon)$ is understood.

\par\medskip\noindent\textbf{Projected equation and projected family.}
For fixed $(\xi,\varepsilon)$, a projected solution is a function $v\in H^1(M)$ such that $v-\chi u_{(\xi,\varepsilon)}\in\E_{(\xi,\varepsilon)}(M)$ and
\[
\langle E_{\widehat g}(v),\psi\rangle=0 \qquad\text{for every }\psi\in\E_{(\xi,\varepsilon)}(M).
\]
We call this equation the projected equation and a family $v_{(\xi,\varepsilon)}$ of such functions a projected family.

\begin{lemma}[{Uniform pairings and projections}]\label{lem:compact-pairings}
The matrix $G:=(\int_M\chi^2\lambda\varphi_{(\xi,\varepsilon,i)} Z_{(\xi,\varepsilon,j)}\,dV_{\widehat g})$ satisfies
\begin{equation}\label{eq:compact-Gram}
\left\|G-\frac{J_\sph}{4(n+1)\bar\varepsilon}\diag(-1,1,\ldots,1)\right\|_{C^1(\overline\Omega,\R^{(n+1)\times(n+1)})} \le C\lambda^n.
\end{equation}
The fields $\chi\varphi_{(\xi,\varepsilon,j)},\chi \partial_{\bar\theta}\varphi_{(\xi,\varepsilon,j)}$ are uniformly bounded in $H^{-1}(M)$, and $\chi\lambda Z_{(\xi,\varepsilon,j)},\chi\lambda \partial_{\bar\theta} Z_{(\xi,\varepsilon,j)}$ are uniformly bounded in $H^1(M)$. In particular the projections
\begin{equation}\label{eq:compact-projection}
f\longmapsto f-\sum_{j,i=0}^n\chi\lambda Z_{(\xi,\varepsilon,j)}(G^{-1})_{ji} \int_M\chi\varphi_{(\xi,\varepsilon,i)} f\,dV_{\widehat g}
\end{equation}
onto $\E_{(\xi,\varepsilon)}(M)$ are uniformly bounded and $C^1$ in normalized parameters.
\end{lemma}
\begin{proof}
Write $\varphi_j:=\varphi_{(\xi,\varepsilon,j)}$ and $Z_j:=Z_{(\xi,\varepsilon,j)}$. On the support of $\chi$, $dV_{\widehat g}=dx$. The full Euclidean pairing is the matrix subtracted from $G$ in \eqref{eq:compact-Gram}, by Lemma~\ref{lem:bubble-pairings}. Outside the fixed ball on which $\chi=1$, $\varphi_0$ satisfies $|\varphi_0|+|\partial_{\bar\theta}\varphi_0|\le C\lambda^{m+2}|y|^{-n-2}$, and the same bound holds for $\varphi_j$, $1\le j\le n$. Also $|\lambda Z_j|+|\partial_{\bar\theta}(\lambda Z_j)|\le C\lambda^m|y|^{2-n}$. Their omitted product has integral $O(\lambda^{2m+2})=O(\lambda^n)$. This proves the pairing estimate including derivatives. Scaling on the inner region and direct tail estimates on the cutoff annulus give the asserted norm bounds. Formula \eqref{eq:compact-projection} and invertibility of $G$ prove the last statement.
\end{proof}

We use the linearization $\LL_{\widehat g,\chi u_{(\xi,\varepsilon)}}$ defined in Section~1. Its bilinear form restricted to $\E_{(\xi,\varepsilon)}(M)$ defines $\LL_{\widehat g,\chi u_{(\xi,\varepsilon)}}^\perp$.

\begin{proposition}[{Uniform projected inverse and parameter estimates}]\label{prop:compact-inverse}
Fix the metric $\widehat g$ and cutoff $\chi$ chosen above. There exist $N_0$ and $C>0$ such that, for $N\ge N_0$ and $(\bar\xi,\bar\varepsilon)\in\overline\Omega$, the operator
\[
\LL_{\widehat g,\chi u_{(\xi,\varepsilon)}}^{\perp}:\E_{(\xi,\varepsilon)}(M) \longrightarrow\E_{(\xi,\varepsilon)}(M)^*
\]
is an isomorphism satisfying, for every $f\in\E_{(\xi,\varepsilon)}(M)$,
\begin{equation}\label{eq:compact-inverse}
\|f\|_{H^1(M)}\le C\sup_{\substack{\psi\in\E_{(\xi,\varepsilon)}(M)\\
\|\psi\|_{H^1(M)}=1}} |\langle\LL_{\widehat g,\chi u_{(\xi,\varepsilon)}}f,\psi\rangle|.
\end{equation}
After identifying nearby subspaces by \eqref{eq:compact-projection}, the inverse is $C^1$ in $(\bar\xi,\bar\varepsilon)$, with its operator norm and first derivative norms bounded by $C$.

More generally, let $G_0\in H^{-1}(M)$ and suppose $f\in H^1(M)$ satisfies
\[
\langle\LL_{\widehat g,\chi u_{(\xi,\varepsilon)}}f,\psi\rangle=\langle G_0,\psi\rangle \quad(\psi\in\E_{(\xi,\varepsilon)}(M)),\qquad \int_M\chi\varphi_{(\xi,\varepsilon,j)}f\,dV_{\widehat g}=d_j.
\]
Then
\begin{align*}
\|f\|_{H^1(M)}&\le C(\|G_0\|_{H^{-1}(M)}+|d|),\\
\|\partial_{\bar\theta}f\|_{H^1(M)}&\le C\bigl(\|\partial_{\bar\theta}G_0\|_{H^{-1}(M)} +\|G_0\|_{H^{-1}(M)}+|\partial_{\bar\theta}d|+|d|\bigr),
\end{align*}
where $d=(d_0,\ldots,d_n)$, $|d|$ is its Euclidean norm. The second estimate applies to $C^1$ data $G_0,d$. All constants are uniform in $N$ and the normalized parameters.
\end{proposition}
\begin{proof}
The uniform projected inverse is proved by contradiction. After rescaling near the concentrating bubble, the projection constraints and Euclidean nondegeneracy force the limiting function to vanish. Local compactness and the decay of the bubble potential then show that its quadratic contribution tends to zero on all of $M$. Testing the linearized equation with the normalized sequence now contradicts coercivity of the conformal Laplacian, which follows from positivity of the Yamabe constant.

Write $u:=u_{(\xi,\varepsilon)}$, $\varphi_j:=\varphi_{(\xi,\varepsilon,j)}$, and $Z_j:=Z_{(\xi,\varepsilon,j)}$. The error estimates at the endpoint orders in Proposition~\ref{prop:bubble-estimates}, cut off on a fixed annulus, give
\begin{equation}\label{eq:compact-approx-Jacobi}
\|E_{\widehat g}(\chi u)\|_{H^{-1}}+ \|\partial_{\bar\theta} E_{\widehat g}(\chi u)\|_{H^{-1}} \le C(\lambda^{l+2}+(\lambda/\rho)^m)=o(1).
\end{equation}
Indeed all coefficients on the cutoff annulus are fixed and bounded, and $u$, its first two spatial derivatives, and its normalized parameter derivatives there are $O(\lambda^m)$. The commutator and nonlinear cutoff terms have $L^q$ norm $O(\lambda^m)\le C(\lambda/\rho)^m$. Differentiating $E_{\widehat g}(\chi u)$ with respect to $\bar\varepsilon,\bar\xi_1,\ldots,\bar\xi_n$ gives $\LL_{\widehat g,\chi u}(\chi\lambda Z_j)=o(1)$ in $H^{-1}$ for $0\le j\le n$.

If \eqref{eq:compact-inverse} fails, there are parameters $(\xi_N,\varepsilon_N)$ and functions $f_N\in\E_{(\xi_N,\varepsilon_N)}(M)$ with $\|f_N\|_{H^1}=1$ and $\|\LL_{\widehat g,\chi u_{(\xi_N,\varepsilon_N)}}^{\perp}f_N \|_{\E_{(\xi_N,\varepsilon_N)}(M)^*}\to0$. Extending that functional to $H^{-1}(M)$ and expressing the difference as a linear combination of the constraint functionals gives
\begin{equation}\label{eq:inverse-contradiction}
\LL_{\widehat g,\chi u}f_N=e_N+\sum_{i=0}^nd_{N,i}\chi\varphi_i, \qquad \|e_N\|_{H^{-1}}\longrightarrow0.
\end{equation}
Testing with $\chi\lambda Z_j$ and using self-adjointness, \eqref{eq:compact-Gram}, and \eqref{eq:compact-approx-Jacobi}, gives $d_{N,i}\to0$ for every $i$.

The rescaled functions $\lambda^m f_N(x_N+\lambda y)$ have bounded $L^2$ gradient and $L^{\frac{2n}{n-2}}$ norms. After passing to a subsequence, the normalized parameters converge, and these functions converge weakly locally in $H^1$ and strongly locally in $L^2$ to $f_\infty\in\dot H^1(\R^n)$. Since $\sup_N\|\varphi_{(\bar\xi_N,\bar\varepsilon_N,j)}\|_{L^q(\R^n\setminus B_R)} \to0$ as $R\to\infty$ for every $j$, local compactness and H\"older's inequality give $f_\infty\in\E_{(\bar\xi_\infty,\bar\varepsilon_\infty)}$. Since the rescaled metric converges to $\delta$ on every fixed ball, passing to the limit in the equation gives $\LL^\perp_{\delta,(\bar\xi_\infty,\bar\varepsilon_\infty)}f_\infty=0$. By the Euclidean isomorphism stated before \eqref{flat-projected-inverse}, this implies $f_\infty=0$. On a ball of radius $R\lambda$, local convergence gives $\int (\chi u)^{p-1}f_N^2\to0$ for fixed $R$. On its complement, H\"older gives
\begin{equation*}
\int_{M\setminus B_{R\lambda}(x_N)}(\chi u)^{p-1}f_N^2\,dV_{\widehat g} \le C\|(\chi u)^{p-1}\mathbf1_{M\setminus B_{R\lambda}(x_N)}\|_{L^{n/2}}.
\end{equation*}
The right side tends to zero as $R\to\infty$, uniformly in large $N$, by the explicit formula for $u$. Thus $\int (\chi u)^{p-1}f_N^2\to0$ on all of $M$. Testing \eqref{eq:inverse-contradiction} with $f_N$ now yields $Q_{\widehat g}(f_N)\to0$, contradicting \eqref{eq:compact-coercivity}. This proves the uniform estimate. For fixed $N$, the restriction of $Q_{\widehat g}$ to $\E(M)$ is coercive and the bounded bubble potential is a compact perturbation from $H^1$ to its dual. The operator $\LL^\perp_{\widehat g,\chi u_{(\xi,\varepsilon)}}$ is therefore Fredholm of index zero. Estimate~\eqref{eq:compact-inverse} gives injectivity and hence surjectivity.

Formula \eqref{eq:compact-projection} identifies nearby subspaces with a fixed subspace, with uniformly bounded first derivatives. The potential $(\chi u)^{p-1}$ is $C^1$ in $L^{n/2}$ and its normalized derivative is uniformly bounded: $|\partial_{\bar\theta}(\chi u)|=|\chi \partial_{\bar\theta}u|\le C\chi u$ and the zero set introduced by $\chi$ is independent of the parameters. The conjugated restricted operator is consequently $C^1$ in operator norm. The inverse-differentiation argument in Lemma~\ref{lem:parameter-derivative-inverse}, applied after the identification \eqref{eq:compact-projection}, gives the derivative estimate.

For prescribed constraints subtract $\sum_{j,i}\chi\lambda Z_j(G^{-1})_{ji}d_i$. This lift and its derivative have $H^1$ norms at most $C(|d|+|\partial_{\bar\theta}d|)$; its image under $\LL_{\widehat g,\chi u}$ obeys the same bound. Apply the homogeneous estimate to the difference.
\end{proof}

\subsection{The projected family and reduced functional}
In this and subsequent formulas, the spatial derivatives of $u_{(\xi,\varepsilon)}$ and $H^\lambda$ are in the translated coordinate $y$.

\begin{lemma}[{Source and cutoff estimates on the manifold}]\label{lem:compact-sources}
For all sufficiently large $N$, the following estimates hold with a constant $C$ independent of $N$. Here $u=u_{(\xi,\varepsilon)}$:
\begin{align}
\|E_{\widehat g}(\chi u)\|_{C^1(\overline\Omega,L^q(M))} &\le C\bigl(\lambda^{l+2}+(\lambda/\rho)^m\bigr),\label{eq:compact-source}\\
\|\chi H^\lambda_{ij}\partial_{ij}u\|_{C^1(\overline\Omega,L^q(M))} &\le C\lambda^{l+2},\label{eq:compact-leading-source}\\
\|E_{\widehat g}(\chi u)-\chi H^\lambda_{ij}\partial_{ij}u\|_{C^1(\overline\Omega,L^q(M))} &\le C\bigl(\lambda^{p(l+2)}+(\lambda/\rho)^m\bigr).\label{eq:compact-improved-source}
\end{align}
The same bound as in \eqref{eq:compact-improved-source} holds for $\LL_{\widehat g,\chi u_{(\xi,\varepsilon)}}(\chi w_{(\xi,\varepsilon)})+\chi H^\lambda_{ij}\partial_{ij}u_{(\xi,\varepsilon)}$. Moreover, for $0\le j\le n$,
\begin{equation}\label{eq:cutoff-w-constraints}
\Bigl\|\int_M\chi^2\varphi_{(\xi,\varepsilon,j)} w_{(\xi,\varepsilon)}\,dV_{\widehat g}\Bigr\|_{C^1(\overline\Omega)} \le C\lambda^n.
\end{equation}
\end{lemma}
\begin{proof}
Write $u:=u_{(\xi,\varepsilon)}$, $w:=w_{(\xi,\varepsilon)}$, and $\varphi_j:=\varphi_{(\xi,\varepsilon,j)}$. For the source in \eqref{eq:compact-leading-source}, scaling gives
\[
\|H^\lambda_{ij}\partial_{ij}u\|_{L^q(\R^n)} =\lambda^{l+2}\|H_{ij} \partial_{ij}u_{(\bar\xi,\bar\varepsilon)}\|_{L^q(\R^n)}.
\]
The $L^q$ norms of $H_{ij}\partial_{ij}u_{(\bar\xi,\bar\varepsilon)}$ and its first parameter derivatives are uniformly bounded: the slowest tail is $O(|y|^{l+2-n})$, which lies in $L^q$ because $(l+2)<m$. Multiplication by the parameter-independent $\chi$ gives \eqref{eq:compact-leading-source}.

On $B_b$ the metric equals $g_0$, so Propositions~\ref{prop:bubble-estimates} and \ref{prop:w-decay}, together with Lemma~\ref{lem:linear-correction-error}, apply with their original constants. Multiplication by $\chi$ cuts off only a fixed region a positive distance from $x_N$. On that region
\[
|\nabla^a\partial_{\bar\theta}^s u|+|\nabla^a\partial_{\bar\theta}^s w|\le C\lambda^m, \qquad  0\le a\le2,\quad s=0,1.
\]
Consequently the commutators with $\chi$, the difference between $\chi u^p$ and $(\chi u)^p$, and the potential-cutoff terms involving $w$ each have $C^1(\overline\Omega,L^q(M))$ norm bounded by $C\lambda^m\le C(\lambda/\rho)^m$. By Lemma~\ref{lem:corrections-full-equations} and scaling, $w$ satisfies \eqref{eq:w-full}. Cutting off that equation gives the asserted estimate. The local metric error is bounded by \eqref{eq:linear-correction-error}.

For the constraints, use $\int_{\R^n}\varphi_jw=0$ and $dV_{\widehat g}=dx$ on the cutoff support to write the integral in \eqref{eq:cutoff-w-constraints} as $\int(\chi^2-1)\varphi_jw$. By Proposition~\ref{prop:w-decay} its absolute value and normalized derivative are bounded by $ C\lambda^{2m+2}\int_{b/2}^\infty r^{l+2-n-1}\,dr\le C\lambda^n. $
\end{proof}

The next lemma and the inequalities following it give the estimates, including first parameter derivatives, used to construct the projected family and compare its normalized reduced functional with $F$.

\begin{lemma}[{Nonlinear remainder estimates}]\label{lem:nonlinear-estimates}
For any $S,\omega,\widetilde\omega\in\E$,
\begin{align*}
\|\NN_S(\omega)\|_{L^{\frac{2n}{n+2}}} &\le C\|\omega\|_\E^{\frac{n+2}{n-2}},\\
\|\NN_S(\omega)-\NN_S(\widetilde\omega)\|_{L^{\frac{2n}{n+2}}} &\le C\bigl(\|\omega\|_\E^{\frac{4}{n-2}}+ \|\widetilde\omega\|_\E^{\frac{4}{n-2}}\bigr) \|\omega-\widetilde\omega\|_\E.
\end{align*}
The constants are independent of $S$. On a fixed closed manifold, the same estimates hold with $\E$ replaced by $H^1(M)$. The map $v\mapsto |v|^{\frac{4}{n-2}}v$ is $C^1$ from $L^{\frac{2n}{n-2}}$ to $L^{\frac{2n}{n+2}}$.
\end{lemma}

\begin{lemma}
Let $u>0$, $v\in\R$, and $n\ge6$. Then
\begin{equation}\label{eq:scalar3}
\left| |v|^{\frac{4}{n-2}}-u^{\frac{4}{n-2}} -\frac{4}{n-2}u^{\frac{6-n}{n-2}}(v-u) \right| \le Cu^{-1}|v-u|^{\frac{n+2}{n-2}}.
\end{equation}
Moreover,
\begin{align}
\left| (|v|^{\frac{4}{n-2}}-u^{\frac{4}{n-2}})uv -\frac{2}n \bigl(|v|^{\frac{2n}{n-2}}-u^{\frac{2n}{n-2}}\bigr) \right| &\le C|v-u|^{\frac{2n}{n-2}},
\label{eq:scalar-energy}\\
\left| \frac{n+2}{n-2}u|v|^{\frac{4}{n-2}} -u^{\frac{n+2}{n-2}} -\frac{4}{n-2}v|v|^{\frac{4}{n-2}} \right| &\le C|v-u|^{\frac{n+2}{n-2}},
\label{eq:scalar-diff1}
\end{align}
and
\begin{equation}\label{eq:scalar-diff2}
\begin{split}
\biggl| &\frac{n-6}{n-2}v|v|^{\frac{4}{n-2}} +\frac{n+2}{n-2}u|v|^{\frac{4}{n-2}} -\frac{n+2}{n-2}u^{\frac{4}{n-2}}v +\frac{6-n}{n-2}u^{\frac{n+2}{n-2}} \biggr|\\
&\hspace{25mm}\le Cu^{-1}|v-u|^{\frac{2n}{n-2}}.
\end{split}
\end{equation}
The constant $C$ depends only on $n$.
\end{lemma}

\begin{proposition}[{The projected family and its linear correction}]\label{prop:compact-projected-family}
There exist $C>0$ and $N_0$ such that, for $N\ge N_0$ and $(\bar\xi,\bar\varepsilon)\in\overline\Omega$, there is a unique $v_{(\xi,\varepsilon)}\in H^1(M)$ satisfying
\begin{equation*}
\begin{gathered}
v_{(\xi,\varepsilon)}-\chi u_{(\xi,\varepsilon)}\in\E_{(\xi,\varepsilon)}(M),\\
\langle E_{\widehat g}(v_{(\xi,\varepsilon)}),\psi\rangle=0 \quad\text{for every }\psi\in\E_{(\xi,\varepsilon)}(M),
\end{gathered}
\end{equation*}
and
\begin{equation}\label{eq:compact-z}
\|v_{(\xi,\varepsilon)}-\chi u_{(\xi,\varepsilon)}\|_{C^1(\overline\Omega,H^1(M))} \le C\bigl(\lambda^{l+2}+(\lambda/\rho)^{\frac{n-2}{2}}\bigr).
\end{equation}
The map $(\bar\xi,\bar\varepsilon)\mapsto v_{(\xi,\varepsilon)}$ is $C^1$ with values in $H^1(M)$. Let $w_{(\xi,\varepsilon)}$ be the linear correction in \eqref{eq:exact-scaling}. Then
\begin{equation}\label{eq:compact-r-bound}
\|v_{(\xi,\varepsilon)}-\chi u_{(\xi,\varepsilon)}-\chi w_{(\xi,\varepsilon)}\|_{C^1(\overline\Omega,H^1(M))} \le C\bigl(\lambda^{\frac{(n+2)(l+2)}{n-2}}+(\lambda/\rho)^{\frac{n-2}{2}}\bigr).
\end{equation}
\end{proposition}
\begin{proof}
Write $u:=u_{(\xi,\varepsilon)}$, $v:=v_{(\xi,\varepsilon)}$, and $w:=w_{(\xi,\varepsilon)}$, and set
\[
z:=v-\chi u,\qquad r:=z-\chi w,\qquad B_2:=\chi H^\lambda_{ij}\partial_{ij}u.
\]
For $z\in\E(M)$, the equation $\langle E_{\widehat g}(v),\psi\rangle=0$ for every $\psi\in\E(M)$ is equivalent to
\[
z=(\LL_{\widehat g,\chi u}^{\perp})^{-1} \bigl(-E_{\widehat g}(\chi u)+\NN_{\chi u}(z)\bigr).
\]
Lemma~\ref{lem:nonlinear-estimates} and the Sobolev inequality give
\[
\|\NN_{\chi u}(z_1)-\NN_{\chi u}(z_2)\|_{H^{-1}} \le C(\|z_1\|_{H^1}+\|z_2\|_{H^1})^{p-1}\|z_1-z_2\|_{H^1}.
\]
Proposition~\ref{prop:compact-inverse} and \eqref{eq:compact-source} therefore make the right-hand side a contraction on a ball of radius $C(\lambda^{l+2}+(\lambda/\rho)^m)$ for large $N$. This proves existence, uniqueness in that ball, and the bound for $\|z\|_{H^1}$. The projections in \eqref{eq:compact-projection} identify nearby subspaces $\E(M)$ with a fixed subspace. The source and inverse depend $C^1$ on the parameters. The nonlinear map is $C^1$ in the parameters and $z$, and its derivative with respect to $z$ has operator norm uniformly less than one, the implicit function theorem gives $C^1$ dependence of $z$, and hence of $v=\chi u+z$, on $(\bar\xi,\bar\varepsilon)$.

To estimate the derivatives,  use
\[
\frac{1}{n(n-2)}\NN_{\chi u}(z)=|\chi u+z|^{p-1}(\chi u+z) -(\chi u)^p-p(\chi u)^{p-1}z,
\]

Where $\chi u>0$,
\[
\begin{aligned}
& \frac{1}{n(n-2)}\partial_{\bar\theta}\NN_{\chi u}(z)=p\bigl(|\chi u+z|^{p-1}-(\chi u)^{p-1}\bigr)\partial_{\bar\theta}z\\
&\quad+p\bigl(|\chi u+z|^{p-1}-(\chi u)^{p-1} -(p-1)(\chi u)^{p-2}z\bigr)\partial_{\bar\theta}(\chi u).
\end{aligned}
\]
Where $\chi u=0$, $\partial_{\bar\theta}(\chi u)=0$ and $ \partial_{\bar\theta}\NN_{\chi u}(z)=n(n-2)p|z|^{p-1}\partial_{\bar\theta}z$. Thus \eqref{eq:scalar3}, $|\partial_{\bar\theta}(\chi u)|\le C\chi u$, and H\"older's inequality give
\begin{equation}\label{eq:compact-nonlinear-derivative}
\|\partial_{\bar\theta} \NN_{\chi u}(z)\|_{L^q}\le C\|z\|_{H^1}^{p-1}\|\partial_{\bar\theta}z\|_{H^1}+C\|z\|_{H^1}^p.
\end{equation}
Apply the derivative estimate of Proposition~\ref{prop:compact-inverse} to $z$, with $G_0=-E_{\widehat g}(\chi u)+\NN_{\chi u}(z)$ and $d_j=0$. Using \eqref{eq:compact-source} and the bound for $\|z\|_{H^1}$ gives $\|\partial_{\bar\theta}z\|_{H^1}\le C(\lambda^{l+2}+(\lambda/\rho)^m) +C\|z\|_{H^1}^{p-1}\|\partial_{\bar\theta}z\|_{H^1}$. Since $\|z\|_{H^1}\to0$ uniformly, absorption proves \eqref{eq:compact-z}. Lemma~\ref{lem:nonlinear-estimates} and \eqref{eq:compact-nonlinear-derivative} then yield
\begin{equation}\label{eq:compact-nonlinear-C1}
\|\NN_{\chi u}(z)\|_{C^1(\overline\Omega,L^q(M))} \le C(\lambda^{p(l+2)}+(\lambda/\rho)^m).
\end{equation}

Since $r=z-\chi w$, we have, for every $\psi\in\E(M)$, that
\[
\langle\LL_{\widehat g,\chi u}r,\psi\rangle =-\langle E_{\widehat g}(\chi u)-B_2,\psi\rangle +\langle\NN_{\chi u}(z),\psi\rangle -\langle\LL_{\widehat g,\chi u}(\chi w)+B_2,\psi\rangle.
\]
Moreover, $\int_M\chi\varphi_{(\xi,\varepsilon,j)}r\,dV_{\widehat g} =-\int_M\chi^2\varphi_{(\xi,\varepsilon,j)}w\,dV_{\widehat g}$ for $0\le j\le n$. Lemma~\ref{lem:compact-sources} and  \eqref{eq:compact-nonlinear-C1} bound the right-hand side in $C^1(\overline\Omega,H^{-1}(M))$ by $C(\lambda^{p(l+2)}+(\lambda/\rho)^m)$; \eqref{eq:cutoff-w-constraints} bounds these integrals in $C^1(\overline\Omega)$ by $C\lambda^n$. Proposition~\ref{prop:compact-inverse} gives $\|r\|_{C^1(\overline\Omega,H^1(M))}\le C(\lambda^{p(l+2)}+(\lambda/\rho)^m+\lambda^n)$. Since $\lambda^n\le(\lambda/\rho)^m$ for large $N$, this proves \eqref{eq:compact-r-bound}.
\end{proof}

\par\medskip\noindent\textbf{The reduced functional.}
For the projected family on the fixed metric $\widehat g$, define the {\it reduced functional} by 
\begin{equation}\label{eq:compact-reduced}
\widehat F_N(\bar\xi,\bar\varepsilon) :=J_{\widehat g}^M(v_{(\xi,\varepsilon)})-J_\sph.
\end{equation}

\begin{proposition}[{$C^1$ comparison for the endpoint orders}]\label{prop:compact-energy-expansion}
Fix an endpoint pair $(n,l)$ in the ranges of Section~\ref{sec:global}. Let $v_{(\xi,\varepsilon)}$ be the family of Proposition~\ref{prop:compact-projected-family}. Recall that $F$ is the limiting functional in \eqref{eq:auxiliary-function}, and define the normalized reduced functional by
\[
\widetilde F_N:=\lambda_N^{-2(l+2)}\widehat F_N.
\]
Then, as $N\to\infty$,
\begin{equation}\label{eq:compact-energy-limit}
\widetilde F_N \longrightarrow F \quad\text{in }C^1(\overline\Omega).
\end{equation}
More precisely,
\begin{equation}\label{eq:compact-energy-error}
\begin{split}
\|\widehat F_N-\lambda^{2(l+2)}F\|_{ C^1(\overline\Omega)} \le C\Bigl[&\lambda^{\frac{2n(l+2)}{n-2}}+\lambda^{l+2}(\lambda/\rho)^{\frac{n-2}{2}}+(\lambda/\rho)^{n-2}\\
&+\lambda^{3(l+2)}(1+\log(\frac{\rho}{\lambda}))^{\mathbf1_{\{n=3l+8\}}}\Bigr].
\end{split}
\end{equation}
The constant $C$ is independent of $N$.
\end{proposition}
\begin{proof}
Write $u:=u_{(\xi,\varepsilon)}$, $v:=v_{(\xi,\varepsilon)}$, and $w:=w_{(\xi,\varepsilon)}$, and set
\[
z:=v-\chi u,\qquad r:=v-\chi u-\chi w,\qquad B_2:=\chi H^\lambda_{ij}\partial_{ij}u.
\]
For $s\ge0$ and $z\in\R$, set
\[
\mathscr C(s,z):=s|s+z|^{p-1}(s+z)-s^p(s+z) -\frac{2}{n}\bigl(|s+z|^{p+1}-s^{p+1}\bigr).
\]
Taking $\psi=z\in\E(M)$ in $\langle E_{\widehat g}(v),\psi\rangle=0$ gives
\[
Q_{\widehat g}(v)-Q_{\widehat g}(\chi u) =\langle E_{\widehat g}(\chi u),z\rangle +n(n-2)\int_M(|v|^{p-1}v+(\chi u)^p)z\,dV_{\widehat g}.
\]
The pointwise identity
\[
(|v|^{p-1}v+(\chi u)^p)(v-\chi u) +( |v|^{p-1}-(\chi u)^{p-1})\chi u\,v =|v|^{p+1}-(\chi u)^{p+1}
\]
and the definition of the Yamabe functional therefore give
\begin{equation}\label{eq:compact-energy-exact}
\begin{aligned}
&J_{\widehat g}^M(v)-J_{\widehat g}^M(\chi u)-\int_M B_2\chi w\,dV_{\widehat g}\\
&\quad=\langle E_{\widehat g}(\chi u)-B_2,z\rangle +\int_M B_2r\,dV_{\widehat g} -n(n-2)\int_M \mathscr C(\chi u,z)\,dV_{\widehat g}.
\end{aligned}
\end{equation}
Lemma~\ref{lem:compact-sources} and Proposition~\ref{prop:compact-projected-family} bound the sum of the $C^1(\overline\Omega)$ norms of the first two terms on the right hand side of \eqref{eq:compact-energy-exact} by

\[
C\Bigl[(\lambda^{p(l+2)}+(\lambda/\rho)^m) (\lambda^{l+2}+(\lambda/\rho)^m) +\lambda^{l+2}(\lambda^{p(l+2)}+(\lambda/\rho)^m)\Bigr].
\]

For the last term, \eqref{eq:scalar-energy} gives $|\mathscr C(\chi u,z)|\le C|z|^{p+1}$. Where $\chi u>0$, direct differentiation gives
\begin{align*}
\mathscr C_z(\chi u,z)&=p\chi u|v|^{p-1}-(\chi u)^p-(p-1)v|v|^{p-1},\\
\mathscr C_s(\chi u,z)&=(2-p)v|v|^{p-1}+p\chi u|v|^{p-1} -p(\chi u)^{p-1}v+(p-2)(\chi u)^p.
\end{align*}
Equations~\eqref{eq:scalar-diff1}--\eqref{eq:scalar-diff2} and $|\partial_{\bar\theta}(\chi u)|=|\chi \partial_{\bar\theta}u|\le C\chi u$ imply
\begin{equation}
|\partial_{\bar\theta}\mathscr C(\chi u,z)|\le C|z|^p|\partial_{\bar\theta}z|+C|z|^{p+1}.
\end{equation}
Where $\chi u=0$, one has $\partial_{\bar\theta}(\chi u)=\chi \partial_{\bar\theta}u=0$ and $\mathscr C(0,z)=-\frac{2}{n}|z|^{p+1}$, so the same bound holds. The $C^1$ nonlinear map in Lemma~\ref{lem:nonlinear-estimates} and H\"older's inequality justify differentiating the integral. Consequently
\[
\left\|\int_M \mathscr C(\chi u,z)\,dV_{\widehat g}\right\|_{C^1(\overline\Omega)} \le C\|z\|_{C^1(\overline\Omega,H^1(M))}^{p+1}.
\]
This proves that the $C^1(\overline\Omega)$ norm of the left side of \eqref{eq:compact-energy-exact} is bounded by $C(\lambda^{(p+1)(l+2)}+\lambda^{l+2}(\lambda/\rho)^m+(\lambda/\rho)^{n-2})$.

The original translated metric $\bar g(y)=g_0(x_N+y)$ is defined on $\R^n$. On the support of $\chi$, the coefficients of $\widehat g$ in the translated coordinates agree with those of $\bar g$. Since all differences below occur a fixed distance from the bubble center, \eqref{eq:bubble-pointwise} gives
\begin{equation}
\|J_{\widehat g}^M(\chi u)-J_{\bar g}^{\R^n}(u)\|_{C^1(\overline\Omega)} \le C\lambda^{n-2}.
\end{equation}
The gradient and scalar terms have this order; the critical nonlinear tail has the smaller order $\lambda^n$. Likewise Proposition~\ref{prop:w-decay} gives
\begin{equation}\label{eq:cutoff-correction-energy}
\left\|\int_M B_2\chi w\,dV_{\widehat g} -\int_{\R^n}H^\lambda_{ij}\partial_{ij}u\,w\,dy\right\|_{C^1(\overline\Omega)} \le C\lambda^{n-2}.
\end{equation}
For \eqref{eq:cutoff-correction-energy}, the slowest radial power in the omitted integral is $r^{2(l+2)+1-n}$, which is integrable at infinity by $n>2(l+2)+2$. A normalized derivative has the same radial majorant.

Estimate~\eqref{eq:polynomial-comparison} gives
\[
\Bigl\|J_{\bar g}^{\R^n}(u)-J_\sph+\int_{\R^n}H^\lambda_{ij}\partial_{ij}u\,w\,dy-\lambda^{2(l+2)}F(\bar\xi,\bar\varepsilon)\Bigr\|_{C^1(\overline\Omega)}=o\bigl(\lambda^{2(l+2)}\bigr).
\]
Combining the quantitative bound \eqref{eq:polynomial-comparison} with \eqref{eq:compact-energy-exact}--\eqref{eq:cutoff-correction-energy}, and using $\lambda^{n-2}\le(\lambda/\rho)^{n-2}$, gives \eqref{eq:compact-energy-error}. Dividing the four displayed error terms by $\lambda^{2(l+2)}$ gives, respectively,
\[
\lambda^{(p-1)(l+2)},\quad \lambda^{m-l-2}\rho^{-m},\quad \lambda^{2(m-l-2)}\rho^{-2m},\quad \lambda^{l+2}(1+\log(\tfrac{\rho}{\lambda}))^{\mathbf1_{\{n=3(l+2)+2\}}}.
\]
All tend to zero since $m>(l+2)$, $\lambda=2^{-N}$, and $\rho=(4N^2)^{-1}$. This proves \eqref{eq:compact-energy-limit}.
\end{proof}

\subsection{Proofs of the main results}
For orders above the endpoint, we localize the whole endpoint tensor $H=T_\loc+T_\glo$, and multiply each rescaled copy by $\lambda_N^{1/4}$. A separate homogeneous tensor of degree $l+2$ prescribes the higher Weyl vanishing order. Thus the localized copies of $H$ determine $F$, while the new homogeneous term contributes only a smaller error.

\begin{proposition}[{Higher prescribed orders}]\label{prop:higher-vanishing-orders}
Fix an endpoint pair $(n,l_0)$, with $H$ and $F$ as in Proposition~\ref{prop:endpoint-critical}. For every integer $l>l_0$, there is a smooth compactly supported small tensor $h_0$ such that $W_{\exp(h_0)}$ vanishes exactly to order $l$ at $0$. For every extension $\widehat g$ supplied by Lemma~\ref{lem:metric-extension} and all sufficiently large $N$, there is a projected family $v_{(\xi,\varepsilon)}$ satisfying
\[
\|v_{(\xi,\varepsilon)}-\chi u_{(\xi,\varepsilon)}\|_{C^1(\overline\Omega,H^1(M))}\le C\lambda_N^{l_0+2+1/4},
\]
where $C$ is independent of $N$.
Its reduced functional $\widehat F_N$, defined by \eqref{eq:compact-reduced}, satisfies, as $N\to\infty$,
\begin{equation}\label{eq:higher-order-energy-limit}
\widetilde F_N:=\lambda_N^{-2(l_0+2)-1/2}\widehat F_N
\longrightarrow F \quad\text{in }C^1(\overline\Omega).
\end{equation}
\end{proposition}
The proof is given in Appendix~\ref{sec:higher-order-proof}.

\begin{lemma}[{Critical points yield positive Yamabe solutions}]\label{lem:compact-critical}
At any critical point of $\widehat F_N$, the corresponding projected solution $v_{(\xi,\varepsilon)}$ is a smooth positive solution of $E_{\widehat g}(v_{(\xi,\varepsilon)})=0$, provided $N$ is sufficiently large. This applies to both Proposition~\ref{prop:compact-energy-expansion} and Proposition~\ref{prop:higher-vanishing-orders}.
\end{lemma}
\begin{proof}
Write $u:=u_{(\xi,\varepsilon)}$, $v:=v_{(\xi,\varepsilon)}$, $\varphi_j:=\varphi_{(\xi,\varepsilon,j)}$, and $Z_j:=Z_{(\xi,\varepsilon,j)}$. Since $E_{\widehat g}(v)$ vanishes on $\E(M)$, there are unique coefficients $d_i$ such that
\[
E_{\widehat g}(v)=\sum_{i=0}^nd_i\chi\varphi_i.
\]
Since $z:=v-\chi u$ satisfies the constraints, differentiation gives, with $\bar\theta=\bar\varepsilon$ for $j=0$ and $\bar\theta=\bar\xi_j$ for $1\le j\le n$,
\begin{equation*}
\int_M\chi\varphi_i\partial_{\bar\theta}v\,dV_{\widehat g} =\int_M\chi^2\lambda\varphi_i Z_j\,dV_{\widehat g} -\int_M(\chi \partial_{\bar\theta}\varphi_i)z\,dV_{\widehat g} =G_{ij}+O(\|z\|_{H^1}).
\end{equation*}
The matrix is invertible for large $N$ by Lemma~\ref{lem:compact-pairings} and either \eqref{eq:compact-z} or \eqref{eq:high-corrections}. The first variation
\[
\partial_{\bar\theta}\widehat F_N=2\sum_{i=0}^nd_i \int_M\chi\varphi_i\partial_{\bar\theta}v\,dV_{\widehat g}
\]
therefore forces every multiplier to vanish at a critical point.

Let $v_-:=\min\{v,0\}$. The exact signed equation and \eqref{eq:compact-coercivity}, tested with $v_-$, imply
\[
c\|v_-\|_{L^{\frac{2n}{n-2}}}^2 \le Q_{\widehat g}(v_-) =n(n-2)\|v_-\|_{L^{\frac{2n}{n-2}}}^{\frac{2n}{n-2}}.
\]
If $v_-$ is nonzero, this inequality imposes a fixed positive lower bound on its critical norm. But $\chi u\ge0$, so $|v_-|\le|z|$ and that norm tends to zero. Thus $v\ge0$. The integral $\int_M|\chi u|^{\frac{2n}{n-2}}\,dV_{\widehat g}$ converges to $\int_{ \R^n} u_{(0,1)}^{ \frac{2n}{n-2} }$ as $N\to\infty$. This, together with $\|z\|_{L^{\frac{2n}{n-2}}(M)}\to0$, also gives $v\not\equiv0$. Trudinger's regularity theorem \cite[Theorem~3]{Trudinger} gives $v\in C^\infty(M)$. Since $v\ge0$ and $v\not\equiv0$, the strong maximum principle on the connected manifold $M$ gives $v>0$.
\end{proof}

\begin{corollary}\label{cor:critical-point-persistence}
For either construction in Propositions~\ref{prop:compact-energy-expansion} and \ref{prop:higher-vanishing-orders}, with endpoint coefficients from Proposition~\ref{prop:endpoint-critical}, $\widehat F_N$ has critical points $(\bar\xi_N,\bar\varepsilon_N)$ for all sufficiently large $N$ such that
\[
(\bar\xi_N,\bar\varepsilon_N)\longrightarrow(0,1),\qquad \widetilde F_N(\bar\xi_N,\bar\varepsilon_N)\longrightarrow F(0,1)<0.
\]
\end{corollary}
\begin{proof}
Choose a ball $B\Subset\Omega$ centered at $(0,1)$ whose closure contains no other critical point of $F$. Nondegeneracy and the $C^1$ convergence in \eqref{eq:compact-energy-limit} and \eqref{eq:higher-order-energy-limit} give, for large $N$,
\[
\deg(\nabla\widetilde F_N,B,0)=\sgn\det D^2F(0,1)\ne0;
\]
see \cite{Nirenberg}. Thus $\widetilde F_N$, and hence its positive scalar multiple $\widehat F_N$, has a critical point $(\bar\xi_N,\bar\varepsilon_N)\in B$. Uniform gradient convergence makes every limit point critical for $F$, hence equal to $(0,1)$. Uniform convergence of $\widetilde F_N$ gives the second limit.
\end{proof}

\begin{proof}[Proof of Theorem~\ref{thm:global-blowup}]
\textbf{Fixing the metric.} Fix the dimension, the prescribed integer $l$, and the point $p\in M$. Choose the metric $g_*$ of Lemma~\ref{lem:positive-flat}. Its chart sends $p$ to $0$. Let $l_0$ be the endpoint order for that dimension. Proposition~\ref{prop:endpoint-critical} supplies the fixed endpoint coefficients and a nondegenerate critical point $(0,1)$ of the limiting functional $F$, with $F(0,1)<0$. A common small rescaling of the two Weyl tensors multiplies $F$ by a positive constant and preserves the critical point and the signature of its Hessian.

If $l=l_0$, use the original accumulating tensor \eqref{eq:h0}, whose smoothness and Weyl vanishing order follow from Lemmas~\ref{lem:smooth-bump-series} and \ref{lem:weyl-vanishing-order}. If $l>l_0$, use \eqref{eq:high-h0}, with the additional fixed background amplitude $\kappa>0$. In either case choose the tensor amplitudes sufficiently small for insertion and then apply Lemma~\ref{lem:metric-extension}. This fixes a smooth metric $\widehat g$ on $M$ with coercive quadratic form $Q_{\widehat g}$, the prescribed Weyl vanishing order at $p$, and the global nonvanishing condition. All these choices are fixed independently of $N$.
\par\medskip\noindent\textbf{Applying the reduced-functional result.}
Corollary~\ref{cor:critical-point-persistence} gives critical points $(\bar\xi_N,\bar\varepsilon_N)\to(0,1)$ of the reduced functionals for this fixed metric. Let $v_N:=v_{(\lambda_N\bar\xi_N,\lambda_N\bar\varepsilon_N)}$ be the member of the projected family on $M$ at these parameters. By Lemma~\ref{lem:compact-critical}, $v_N$ is a smooth positive solution of \eqref{eq:yamabe}.
\par\medskip\noindent\textbf{Blow-up at $p$.}
Its center is $x_N+\lambda_N\bar\xi_N\to0$ and its scale is $\lambda_N\bar\varepsilon_N$, where $\bar\varepsilon_N\to1$. Fix $R>0$. On $B_{R\lambda_N}(x_N+\lambda_N\bar\xi_N)$, the cutoff is one, and \eqref{eq:compact-z} or \eqref{eq:high-corrections}, together with Sobolev, gives
\[
\|v_N-u_{(\lambda_N\bar\xi_N,\lambda_N\bar\varepsilon_N)}\|_{L^{\frac{2n}{n-2}}}\to0.
\]
By scaling and $\bar\varepsilon_N\to1$, the $L^{\frac{2n}{n-2}}$ norm of the bubble on this ball is bounded below by a positive constant, and hence so is that of $v_N$. Since the ball has volume $O(\lambda_N^n)$, Hölder's inequality yields
\[
\sup_{B_{R\lambda_N}(x_N+\lambda_N\bar\xi_N)}v_N \ge c\lambda_N^{-\frac{n-2}{2}}\to\infty.
\]

\par\medskip\noindent\textbf{The volume conclusion.}
The Yamabe functional at the critical point has the expansion
\begin{equation}\label{eq:critical-functional-values}
J_{\widehat g}^M(v_N)-J_\sph
=\begin{cases}
\lambda_N^{2(l_0+2)}\bigl(F(0,1)+o(1)\bigr),& l=l_0,\\
\lambda_N^{2(l_0+2)+1/2}\bigl(F(0,1)+o(1)\bigr),& l>l_0.
\end{cases}
\end{equation}
It is negative for sufficiently large $N$ and tends to zero. Set $g_N:=v_N^{\frac{4}{n-2}}\widehat g$. By \eqref{eq:conformal-law}, these metrics have scalar curvature $4n(n-1)$. Testing the equation for $v_N$ with $v_N$ gives
\[
Q_{\widehat g}(v_N)=n(n-2)\int_Mv_N^{\frac{2n}{n-2}}\,dV_{\widehat g}.
\]
Since the integral is $\operatorname{Vol}(M,g_N)$, the definitions yield
\begin{equation}\label{eq:compact-energy-levels}
J_{\widehat g}^M(v_N)=2(n-2)\operatorname{Vol}(M,g_N).
\end{equation}
The expansion \eqref{eq:critical-functional-values} gives $J_{\widehat g}^M(v_N)<J_\sph$ and $J_{\widehat g}^M(v_N)\to J_\sph$. By \eqref{eq:spherical-level} and \eqref{eq:compact-energy-levels}, this is equivalent to $\operatorname{Vol}(M,g_N)<\frac{J_\sph}{2(n-2)}$ and $\operatorname{Vol}(M,g_N)\to\frac{J_\sph}{2(n-2)}$. The volumes therefore lie below and converge to the volume of the round sphere with the same scalar curvature. The exact Weyl order follows from Lemma~\ref{lem:weyl-vanishing-order} or \eqref{eq:high-jet-exact}, since the final metric is unchanged on the ball $B_{3b}$. The global Weyl condition was imposed in Lemma~\ref{lem:metric-extension}. Discard finitely many initial indices and relabel by $\nu$.
\end{proof}

\appendix
\section{Curvature and integral identities}\label{sec:identities}
We record the curvature expansion and integral identities used in the Euclidean estimates and the calculations of the limiting functional, as well as some basic facts about the Aubin-Talenti bubbles.
\subsection{Identities for Aubin--Talenti bubbles}
\label{subsec:bubble-identities}

Recall that for $(\xi,\varepsilon)\in\R^n\times(0,\infty)$, we let $ u_{(\xi,\varepsilon)}(x)$ be as in (\ref{TheBubble}).

The following integral identities and norm estimates follow by direct computation at $(\xi,\varepsilon)=(0,1)$, followed by translation and scaling.

\begin{lemma}\label{lem:bubble-pairings}
The following identities hold:
\begin{align}
\left(\int_{\R^n}\varphi_{(\xi,\varepsilon,i)}Z_{(\xi,\varepsilon,j)}\dd\right)_{i,j=0}^n &=\frac{J_\sph}{4(n+1)\varepsilon}\diag(-1,1,\ldots,1),\label{eq:dual-Jacobi}\\
\left(\int_{\R^n}\varphi_{(\xi,\varepsilon,i)}\varphi_{(\xi,\varepsilon,j)}\dd\right)_{i,j=0}^n &=\varepsilon^{-2}\diag(\|\varphi_{(0,1,0)}\|_2^2,\ldots, \|\varphi_{(0,1,n)}\|_2^2).
\label{eq:dual-dual}
\end{align}
For every parameter $\theta\in\{\varepsilon,\xi_1,\ldots,\xi_n\}$, $\|\varphi_{(\xi,\varepsilon,j)}\|_2\ge C^{-1}\varepsilon^{-1}$ and

\begin{equation}\label{eq:field-norms}
\begin{aligned}
\|\varphi_{(\xi,\varepsilon,j)}\|_{L^{\frac{2n}{n+2}}}+\varepsilon\|\partial_\theta\varphi_{(\xi,\varepsilon,j)}\|_{L^{\frac{2n}{n+2}}} &\le C,\\
\varepsilon \|\varphi_{(\xi,\varepsilon,j)}\|_\E+ \varepsilon^2\|\partial_\theta\varphi_{(\xi,\varepsilon,j)}\|_\E + \|Z_{(\xi,\varepsilon,j)}\|_\E+\|Z_{(\xi,\varepsilon,j)}\|_{L^{\frac{2n}{n-2}}}&\le C\varepsilon^{-1}.
\end{aligned}
\end{equation}

\end{lemma}

\subsection{General identities}

The following lemma is from Ambrosetti--Malchiodi \cite[Lemma~2.1]{AmbrosettiMalchiodi} and Brendle \cite[Proposition~26]{Brendle}.
\begin{lemma}\label{lem:scalar-expansion}
Let $h$ be symmetric and trace-free, with $|h|\le1$, and let $g=\exp(h)$, where $\exp$ is the matrix exponential. Then
\begin{equation*}
\begin{split}
\left|R_g-\partial_i\partial_jh_{ij} +\partial_i(h_{ik}\partial_jh_{jk}) -\frac{1}{2}\partial_ih_{ik}\partial_jh_{jk} +\frac{1}{4}\partial_kh_{ij}\partial_kh_{ij}\right|\\
\le C\bigl(|h|^2|\partial^2h|+|h||\partial h|^2\bigr).
\end{split}
\end{equation*}
In a region where $\partial_ih_{ij}=0$,
\begin{equation*}
R_g=-\frac{1}{4}|\partial h|^2 +O(|h|^2|\partial^2h|+|h||\partial h|^2).
\end{equation*}
\end{lemma}
\begin{lemma}
If $P,Q$ are homogeneous polynomials of degrees $k,m$ and $Q$ is harmonic, then
\[
\int_{S^{n-1}}\nabla P\cdot\nabla Q =m(m+k+n-2)\int_{S^{n-1}}PQ.
\]
If $F_i$ are homogeneous of degree $k$, then
\[
\int_{S^{n-1}}x_iF_i =\int_{B_1}\partial_iF_i =\frac{1}{n+k-1}\int_{S^{n-1}}\partial_iF_i.
\]
In particular, for homogeneous $P$ of degree $k\ge2$,
\[
\int_{S^{n-1}}P=\frac{1}{k(n+k-2)}\int_{S^{n-1}}\Delta P.
\]
\end{lemma}
We write $\mathrm B(a,b):=\Gamma(a)\Gamma(b)/\Gamma(a+b)$ for the beta function, where $\Gamma$ is Euler's gamma function. The beta-integral recurrences in \cite[Proposition~27]{Brendle} imply, whenever both sides converge,
\begin{equation*}
\begin{gathered}
\mci_j^i=\frac{2j-i-3}{2j-2}\mci_{j-1}^i =\frac{2j-i-3}{i+1}\mci_j^{i+2} =\frac{i-1}{2j-i-1}\mci_j^{i-2} =\frac{i-1}{2j-2}\mci_{j-1}^{i-2},\\
\int_0^\infty\frac{(1-r^2)r^i}{(1+r^2)^j}dr =\frac{2j-2i-4}{2j-i-3}\mci_j^i.
\end{gathered}
\end{equation*}

We use the spherical moment identities of degrees two, four, and six recorded in \cite[Corollary~29]{Brendle}.

If $F(x')$ is homogeneous of degree $k$, then
\begin{equation*}
\int_{S^{n-1}}F(x')=2(n+k)\mci_{(n+k+2)/2}^{n+k-2} \int_{S^{n-2}}F =\frac{n+k}{n+k-1} \int_{S^{n-1}}F(x')|x'|^2.
\end{equation*}

\begin{lemma}\label{lem:weyl-contractions}
Let $W\in\Wspace$ satisfy $W_{1abc}=0$. Then
\begingroup\allowdisplaybreaks[4]
\begin{equation*}
\int_{S^{n-1}}\Phi(W)_{ij}^2=\mu_n|W|^2,\qquad \int_{S^{n-1}}(\partial_k\Phi(W)_{ij})^2=2(n+2)\mu_n|W|^2.
\end{equation*}
\begin{align}
\int_{S^{n-1}}\Phi(W)_{pk}\Phi(W)_{qk}&=\mu_n(W\circ W)_{pq},\label{eq:QQ-rank2}\\
\int_{S^{n-1}}(\partial_k\Phi(W)_{ij})^2x_px_q &=4\mu_n(W\circ W)_{pq}+2\mu_n|W|^2\delta_{pq},\label{eq:dQ-rank2}\\
\int_{S^{n-1}}\Phi(W)_{ij}^2x_px_q &=\frac{\mu_n}{n+4}\bigl(4(W\circ W)_{pq}+|W|^2\delta_{pq}\bigr).
\label{eq:Q-rank2}
\end{align}
\endgroup
For $V\in\Wspace$ with $V_{1abc}=0$, polarization yields
\begin{equation*}
\int_{S^{n-1}}\ip{\Phi(W)}{\Phi(V)}=\mu_n\ip WV.
\end{equation*}
Thus orthogonal Weyl tensors have quadratic images orthogonal in $L^2(S^{n-1})$. The corresponding derivative contraction also vanishes.
\end{lemma}
For nonnegative integers $A,B$ and $r\in\{1,2,3\}$, define
\begin{equation}\label{eq:Rratio}
R_r(A,B):=\mathrm B(A+\frac{1}{2},B+\frac{n-1}{2}+r)/\mathrm B(\frac{1}{2},\frac{n-1}{2}+r).
\end{equation}
If either argument is negative, set $R_r(A,B)=0$.

\begin{lemma}\label{lem:weighted-Q}
For $W\in\Wspace$ with $W_{1abc}=0$ and nonnegative integers $A,B$,
\begin{align*}
\int_{S^{n-1}}t^{2A}\varrho^{2B}\Phi(W)_{ij}^2 &=\mu_n|W|^2R_2(A,B),\\
\int_{S^{n-1}}t^{2A}\varrho^{2B}(\partial_k\Phi(W)_{ij})^2 &=2(n+2)\mu_n|W|^2R_1(A,B).
\end{align*}
\end{lemma}
\begin{proof}
Put $Q:=\Phi(W)$ and write $x=(t,\varrho\omega)$, where $\varrho=(1-t^2)^{1/2}$ and $\omega\in S^{n-2}$. Then $d\sigma=\varrho^{n-3}dt\,d\sigma(\omega)$, $Q(\varrho\omega)=\varrho^2Q(\omega)$, and $\partial_\alpha Q(\varrho\omega)=\varrho\partial_\alpha Q(\omega)$ for $\alpha\ge2$, while $\partial_1Q=0$. The $t$ integrals are the beta integrals in \eqref{eq:Rratio}. Divide by the unweighted identities in Lemma~\ref{lem:weyl-contractions} to obtain the two ratios.
\end{proof}

\subsection{ Weighted contractions}
We compute the weighted contractions and radial ratios used in Section~\ref{sec:auxiliary}.

\begin{proposition}\label{prop:angular-contractions}
The contractions in \eqref{eq:Q-sum}--\eqref{eq:M-sum} satisfy
\begin{multline}\label{eq:I-angular}
I(a,b,c,d)=\mu_n|W|^2\Bigl(2(n+2)R_1(a+c,b+d)\\
+4acR_2(a+c-1,b+d)+4(bd+b+d)R_2(a+c,b+d-1)\Bigr),
\end{multline}
\begin{equation}\label{eq:A-angular}
A_{pq}(a,b,c,d)=\mu_n(W\circ W)_{pq}R_2(a+c,b+d),
\end{equation}
\begingroup\allowdisplaybreaks[4]
\begin{multline}
B_{pq}(a,b,c,d)\\*
\shoveleft{\quad=\mu_n|W|^2\Bigl[4acR_2(a+c,b+d)+4(bd+b+d)R_2(a+c+1,b+d-1)}\\*
\shoveleft{\hspace{32mm}+2(n+2)R_1(a+c+1,b+d)\Bigr]\delta_{p1}\delta_{q1}}\\[3pt]
\shoveleft{\qquad+\mu_n|W|^2\biggl[2R_2(a+c,b+d)+\frac{4}{n+4}\bigl(acR_3(a+c-1,b+d)}\\*
\shoveleft{\hspace{35mm}+(bd+b+d)R_3(a+c,b+d-1)\bigr)\biggr](\delta_{pq}-\delta_{p1}\delta_{q1})}\\[3pt]
\shoveleft{\qquad+\mu_n\biggl[4R_2(a+c,b+d)+\frac{16}{n+4}\bigl(acR_3(a+c-1,b+d)}\\*
\shoveleft{\hspace{35mm}+(bd+b+d)R_3(a+c,b+d-1)\bigr)\biggr](W\circ W)_{pq}.}\label{eq:B-angular}
\end{multline}
\endgroup
\end{proposition}
\begin{proof}
In this proof put $A:=a+c$, $B:=b+d$, and $Q:=\Phi(W)$.

The derivatives of $H^{a,b}=t^{2a}\varrho^{2b}Q$ are
\[
\partial_1H^{a,b}=2at^{2a-1}\varrho^{2b}Q, \quad \partial_\alpha H^{a,b}=2bt^{2a}\varrho^{2b-2}x_\alpha Q +t^{2a}\varrho^{2b}\partial_\alpha Q.
\]
Use $x'\cdot\nabla Q=2Q$ to obtain the pointwise identity
\[
\begin{aligned}
\partial_kH^{a,b}_{ij}\partial_kH^{c,d}_{ij} ={}&4ac\,t^{2A-2}\varrho^{2B}Q_{ij}^2 +4 (bd+b+d)\,t^{2A}\varrho^{2B-2}Q_{ij}^2\\
&+t^{2A}\varrho^{2B}(\partial_kQ_{ij})^2.
\end{aligned}
\]
Lemma~\ref{lem:weighted-Q} gives \eqref{eq:I-angular}. Applying the same beta integral to \eqref{eq:QQ-rank2} gives \eqref{eq:A-angular}.

Multiply the pointwise identity by $x_px_q$. For $p=q=1$, Lemma~\ref{lem:weighted-Q} gives  the term proportional to $\delta_{p1}\delta_{q1}$ in \eqref{eq:B-angular}. For $p,q>1$, the weighted versions of \eqref{eq:dQ-rank2} and \eqref{eq:Q-rank2} are
\begin{align*}
\int t^{2A}\varrho^{2B}x_px_qQ_{ij}^2 &=\frac{\mu_n}{n+4} \bigl(|W|^2\delta_{pq}+4(W\circ W)_{pq}\bigr)R_3(A,B),\\
\int t^{2A}\varrho^{2B}x_px_q(\partial_kQ_{ij})^2 &=\mu_n\bigl(2|W|^2\delta_{pq}+4(W\circ W)_{pq}\bigr)R_2(A,B).
\end{align*}
They follow by slicing at fixed $t$ and normalizing by the unweighted identities. Substitution gives the coefficients of $\delta_{pq}-\delta_{p1}\delta_{q1}$ and $(W\circ W)_{pq}$ in \eqref{eq:B-angular}. The remaining mixed components vanish by primed reflection.
\end{proof}

\begin{proposition}\label{prop:odd-tensor-contribution}
For $H_V:=\beta t\Phi(V)$, define $I_V:=\int_{S^{n-1}}|\partial H_V|^2$ and
\[
(A_V)_{pq}:=\int_{S^{n-1}}(H_V)_{pk}(H_V)_{qk},\qquad (B_V)_{pq}:=\int_{S^{n-1}}|\partial H_V|^2x_px_q.
\]
Then $I_V=3\mu_n\beta^2|V|^2$ and
$$
\begin{array}{lr}
(A_V)_{11}=0, \ \ \qquad\quad \ \ (A_V)_{pq}=\frac{\mu_n\beta^2}{n+4}(V\circ V)_{pq}& (p,q>1),\\
(B_V)_{11}=\frac{7\mu_n\beta^2}{n+4}|V|^2, \ \ (B_V)_{pq}=\frac{\mu_n\beta^2}{n+4} \bigl(8(V\circ V)_{pq}+3|V|^2\delta_{pq}\bigr)&(p,q>1).
\end{array}
$$
Consequently
\begin{equation}\label{eq:F-V}
F(0,\varepsilon)=-\frac{c(n)}{4}Q_l(\tau,\varepsilon) -\frac{3(n-2)}{16(n-1)}\mu_n|V|^2\beta^2 \mci_{n-2}^{n+3}\varepsilon^6.
\end{equation}
\end{proposition}
\begin{proof}
Writing $Q:=\Phi(V)$, one has $|\partial(tQ)|^2=|Q|^2+t^2|\partial Q|^2$. Lemma~\ref{lem:weighted-Q} gives $I_V=\mu_n\beta^2|V|^2(1+2)=3\mu_n\beta^2|V|^2$. Equations~\eqref{eq:QQ-rank2}--\eqref{eq:Q-rank2}, with weights evaluated using \eqref{eq:Rratio}, give $A_V$ and the transverse part of $B_V$. In the distinguished direction,
\[
\begin{split}
(B_V)_{11} &=\beta^2\int_{S^{n-1}}t^2(|Q|^2+t^2|\partial Q|^2)\\
&=\mu_n\beta^2|V|^2 \bigl(R_2(1,0)+2(n+2)R_1(2,0)\bigr) =7\mu_n\beta^2|V|^2/(n+4).
\end{split}
\]
The $W$--$V$ cross term in the scalar contribution vanishes by parity in $t$. Separate radial and angular variables in \eqref{eq:F-axis-value} to obtain \eqref{eq:F-V}.
\end{proof}

\section{Coefficient formulas and signs}\label{sec:coefficients}
We record all coefficient lists used in Section~\ref{sec:auxiliary}. They are obtained by substituting  the three choices of $H$ in \eqref{eq:Q-sum} and \eqref{eq:M-sum}, using Proposition~\ref{prop:angular-contractions} and the beta-integral recurrences in Appendix~\ref{sec:identities}. Lemma~\ref{lem:endpoint-coefficient-signs} supplies the signs used in Proposition~\ref{prop:endpoint-critical}: positivity of the transverse and scale blocks, and of the distinguished block when $l=6$. The coefficient lists below support that sign calculation.

\begingroup
\allowdisplaybreaks[4]
\setlength{\jot}{2pt}
For $l=2$, the scalar coefficients are
$$
\qhat_0^2=2(n+2)\tau^2,\quad \qhat_1^2=\frac{4(n+2)(n+3)}{n-8}\tau,\quad \qhat_2^2=\frac{2(n+2)(n+3)(n+7)}{(n-10)(n-8)}.
$$
The translation coefficients are
\begin{align}
X_2(\tau)&=-16\bigl(2(n-8)\tau+3(n+5)\bigr)/[(n+4)(n-8)],\label{eq:XYZ2-table}\\
Y_2(\tau)&=-8\bigl((n-8)\tau+n+6\bigr)/[(n+4)(n-8)],\qquad Z_2(\tau)=0.\notag
\end{align}
For $l=4$, the scalar coefficients are
\begin{align}
\qhat_0^4&=2(n+2)\tau^2,\qquad \qhat_1^4=-28(n+2)(n+3)\tau /(n-8),\label{eq:q4-table}\\
\qhat_2^4&=2(n+2)(n+3)\bigl(2(n+5)\tau+49(n+7)\bigr)/[(n-10)(n-8)],\notag\\
\qhat_3^4&=-28(n+2)(n+3)(n+5)(n+11)/[(n-12)(n-10)(n-8)],\notag\\
\qhat_4^4&=\frac{2(n+2)(n+3)(n+5)(n+7)(n+17)}{(n-14)(n-12)(n-10)(n-8)}.\notag
\end{align}
The translation coefficients are
\begin{align}
X_4(\tau)&=\frac{16}{(n+4)(n-8)(n-10)(n-12)}\Bigl[(10n^3-352n^2+4104n\label{eq:XYZ4-table}\\*
&\qquad-15840)\tau-85n^3+2331n^2-10244n-120120\Bigr],\notag\\
Y_4(\tau)&=\frac{8}{(n+4)(n-8)(n-10)(n-12)}\Bigl[(5n^3-176n^2+2052n\notag\\*
&\qquad-7920)\tau-30n^3+786n^2-2523n-48510\Bigr],\notag\\
Z_4(\tau)&=0.\notag
\end{align}
For $l=6$, the scalar coefficients are
\begin{align}
\qhat_0^6&=2(n+2)\tau^2,\qquad \qhat_1^6=4(n+2)(95n+333)\tau/(n-8),\label{eq:q6-table}\\
\qhat_2^6&=\frac{2(n+2)\bigl[9025n^2+99370n+228405-(38n^2+328n+672)\tau\bigr]}{(n-10)(n-8)},\notag\\
\qhat_3^6&=\frac{4(n+2)}{(n-12)(n-10)(n-8)} \Bigl[(n^3+15n^2+71n+105)\tau\notag\\*
&\qquad-1805n^3-36347n^2-208916n-354840\Bigr],\notag\\
\qhat_4^6&=\frac{2(n+2)}{(n-14)(n-12)(n-10)(n-8)} \Bigl[551n^4+17804n^3\notag\\*
&\qquad+185320n^2+764140n+1066242\Bigr],\notag\\
\qhat_5^6&=-\frac{4(n+2)(n+3)(n+5)(n+7)(19n^2+620n+4128)}{(n-16)(n-14)(n-12)(n-10)(n-8)},\notag\\
\qhat_6^6&=\frac{2(n+2)(n+3)(n+5)(n+7)(n+9)(n+11)(n+31)}{(n-18)(n-16)(n-14)(n-12)(n-10)(n-8)}.\notag
\end{align}

The translation coefficients are
\begin{align}
X_6(\tau)&=-\frac{16}{(n+4)(n-8)(n-10)(n-12)(n-14)(n-16)}\Bigl[\bigl(120n^5\notag\\*
&\qquad-8032n^4+215194n^3-2886788n^2+19401296n-52281600\bigr)\tau\label{eq:XYZ6-table}\\
&\qquad+12840n^5-772672n^4+16809106n^3\notag\\*
&\qquad-134347880n^2-207140197n+6173139483\Bigr],\notag\\
Y_6(\tau)&=-\frac{8}{(n+4)(n-8)(n-10)(n-12)(n-14)(n-16)}\Bigl[\bigl(60n^5\notag\\*
&\qquad-4016n^4+107597n^3-1443394n^2+9700648n-26140800\bigr)\tau\notag\\
&\qquad+4620n^5-272032n^4+5666485n^3\notag\\*
&\qquad-39465782n^2-161299159n+2500607508\Bigr],\notag\\
Z_6(\tau)&=-\frac{48}{(n+4)(n-8)(n-10)(n-12)(n-14)}\notag\\*
&\qquad\Bigl[\bigl(6n^4-291n^3+5268n^2-42180n+126000\bigr)\tau\notag\\*
&\qquad+462n^4-18687n^3+222738n^2-128539n-8755495\Bigr].\notag
\end{align}

\endgroup

\begin{lemma}[{Signs determining the endpoint Hessian}]\label{lem:endpoint-coefficient-signs}
Let $\tau_{l,-}<\tau_{l,+}$ be the roots of \eqref{eq:P2}--\eqref{eq:P6} in the corresponding dimension ranges. At the roots $\tau_{2,-}$, $\tau_{4,+}$, and $\tau_{6,-}$ one has
\[
-\sum_{s=0}^l(2s+4)(2s-2)\qhat_s^l(\tau)>0,\qquad Y_l(\tau)+\min(X_l(\tau),0)>0.
\]
Moreover, $Z_2=Z_4=0$, and $Z_6(\tau_{6,-})>0$. The polynomial $P_2$ is negative immediately to the right of $\tau_{2,-}$, and $P_4$ is negative immediately to the left of $\tau_{4,+}$.
\end{lemma}

\begin{proof}
For $l=2$, the discriminant is $(n-10)(n+3)(n^2-55n+178)>0$ for $n\ge52$. Write the roots as $-r_+<-r_-<0$. Then $(n-8)r_+>\frac{3}{2}(n+3)$ and $-\sum_{s=0}^2(2s+4)(2s-2)\qhat_s^2(-r_+)=16(n+2)r_+(r_+-r_-)>0$. Substitution in \eqref{eq:XYZ2-table} gives $Y_2(-r_+)>0$ and $X_2(-r_+)+Y_2(-r_+)>0$, hence $Y_2(-r_+)+\min(X_2(-r_+),0)>0$.

For the remaining finite dimension ranges, direct substitution in \eqref{eq:P4}--\eqref{eq:P6} gives $P_4(21/2)<0<P_4(23/2)$ and $P_6(-143)>0>P_6(-133)$, respectively. Thus $ 21/2<\tau_{4,+}<23/2, -143<\tau_{6,-}<-133. $

For every integer $30\le n\le51$, the formulas \eqref{eq:q4-table} and \eqref{eq:XYZ4-table} give positive values and positive $\tau$-derivatives of $X_4$, $Y_4$, and $-\sum_{s=0}^4(2s+4)(2s-2)\qhat_s^4$ at $\tau=21/2$. For every integer $25\le n\le29$, the formulas \eqref{eq:q6-table} and \eqref{eq:XYZ6-table} give positive values and negative $\tau$-derivatives of $X_6$, $Y_6$, $Z_6$, and $-\sum_{s=0}^6(2s+4)(2s-2)\qhat_s^6$ at $\tau=-133$. These are finite checks with rational numbers. The translation coefficients are affine in $\tau$, and the coefficient of $\tau^2$ in each scale coefficient $-\sum_{s=0}^l(2s+4)(2s-2)\qhat_s^l$ is $16(n+2)>0$. The signs therefore persist from these endpoints to the respective roots. Finally, $Z_2=Z_4=0$ follows from the tables, and the positive leading coefficients and distinct roots of $P_2,P_4$ give the asserted one-sided negativity.
\end{proof}

\section{Threshold minors in the vanishing estimates}\label{sec:KMS}
This appendix records the dimension thresholds underlying the known vanishing statements in Proposition~\ref{prop:fixed-order-vanishing} and the introduction's table. In the coordinate expansion used in \cite[Appendix A]{KMS}, retaining metric degrees $2$ through $5$ corresponds to two radial modes in the degree-two divergence-free sector. Retaining degrees through $7$ corresponds to three such modes. Let $b_0>0$ be the coefficient in \cite[Lemma A.1]{KMS}. The coefficients used in the proof of \cite[Proposition A.4]{KMS} are
\begin{equation}\label{eq:KMS-coeff}
b_{2m}=b_0\prod_{j=1}^m\frac{n+2j-4}{n-2j}, \qquad c_{2m}=\frac{4m}{n-2m-2}b_{2m},
\end{equation}
Define
\begin{equation}\label{eq:KMS-matrix}
M_2^{(3)}:=\frac12
\begin{pmatrix}
(n+2)c_4&(n+4)c_6&(n+6)c_8\\
(n+4)c_6&(n+8)c_8&(n+12)c_{10}\\
(n+6)c_8&(n+12)c_{10}&(n+18)c_{12}
\end{pmatrix}.
\end{equation}
Let $M_2^{(2)}$ be its leading $2\times2$ principal block, and let $\Delta_j$ be the leading $j\times j$ principal minor of $M_2^{(3)}$, $j=1,2,3$.

\begin{lemma}
For $n\ge25$, one has $\Delta_1>0$ and
\begin{align}
\Delta_2&=C_2(n)(-n^2+54n-152),\label{eq:KMS-D2}\\
\Delta_3&=C_3(n)(-n^2+34n-144),\label{eq:KMS-D3}
\end{align}
where
\begin{equation}\label{eq:KMS-C2}
C_2(n):=\frac{4n^2(n+2)^2(n+4)b_0^2}{(n-4)^2(n-6)^2(n-8)^2(n-10)}>0
\end{equation}
and
\begin{equation}\label{eq:KMS-C3}
C_3(n):=\frac{6144n^3(n+2)^3(n+4)^2(n+6)(n-2)(n-3)b_0^3}{(n-4)^3(n-6)^3(n-8)^3(n-10)^3(n-12)^2(n-14)}>0.
\end{equation}
Consequently $M_2^{(2)}$ is positive definite for $25\le n\le51$, and $M_2^{(3)}$ is positive definite for $25\le n\le29$.
\end{lemma}
\begin{proof}
Substitute \eqref{eq:KMS-coeff} into \eqref{eq:KMS-matrix}, clear the positive denominators, and expand the leading $2\times2$ and $3\times3$ determinants. This gives \eqref{eq:KMS-D2}--\eqref{eq:KMS-C3}; it also gives $\Delta_1=(n+2)c_4/2>0$. The polynomial factor in \eqref{eq:KMS-D2} is positive on the stated integer range and negative at $n=52$:
\[
-51^2+54\cdot51-152=1, \qquad -52^2+54\cdot52-152=-48.
\]
The polynomial is concave, so its minimum on $[25,51]$ occurs at an endpoint; both endpoint values are positive. Similarly,
\[
-29^2+34\cdot29-144=1, \qquad -30^2+34\cdot30-144=-24,
\]
and the polynomial factor in \eqref{eq:KMS-D3} is positive on $[25,29]$. The positive factors in \eqref{eq:KMS-C2}--\eqref{eq:KMS-C3} and Sylvester's criterion prove the assertions.
\end{proof}

The Weyl derivatives through order $j$ depend on metric derivatives through order $j+2$. The truncations through degrees $5$ and $7$ in the cited argument give the vanishing conclusions through orders $3$ and $5$.

\section{Higher prescribed orders of vanishing}\label{sec:higher-order-proof}

\begin{proof}[Proof of Proposition~\ref{prop:higher-vanishing-orders}]
\textbf{Metric and Weyl vanishing order.} We construct a metric with the prescribed first nonzero Weyl derivative.

Use the parameters and parameter derivatives of Section~\ref{sec:global}, with $m,p,q$ as fixed in Section~\ref{sec:euclidean-estimates}. Let $T_\glo^{[l]}$ have only the following nonzero components:
\begin{equation*}
(T_\glo^{[l]})_{33}=x_2^{l+2},\qquad (T_\glo^{[l]})_{44}=-x_2^{l+2}.
\end{equation*}
It is trace-free, divergence-free, and independent of $x_1$. The linearized curvature formula \eqref{eq:linear-curvature} gives
\begin{equation*}
(DW_\delta[T_\glo^{[l]}])_{2323} =-\frac{n-3}{2(n-2)}(l+2)(l+1)x_2^l.
\end{equation*}

Keep the original $x_N,\lambda_N,\rho_N,\eta$ and the full endpoint tensor $H$. Set $\mu_N=\lambda_N^{1/4}$, $a_N=\mu_N\lambda_N^{l_0+2}$, and, for a fixed small $\kappa>0$, define
\begin{equation}\label{eq:high-h0}
h_0(x)=\kappa\eta(|x|^2)T_\glo^{[l]}(x)+\sum_{N=N_0}^{\infty}\eta(N^4|x-x_N|^2)a_NH(\lambda_N^{-1}(x-x_N)).
\end{equation}
For every homogeneous degree $j\le l_0+2$, the coefficient in the sum is $2^{-N(l_0+2-j+1/4)}$. The argument of Lemma~\ref{lem:smooth-bump-series} therefore makes the sum smooth and flat at $0$. As in Lemma~\ref{lem:weyl-vanishing-order}, the nonlinear curvature terms have degree at least $2l+2>l$, so
\begin{equation}\label{eq:high-jet-exact}
\begin{gathered}
\nabla^jW_{\exp(h_0)}(0)=0\quad(j<l),\\
(\nabla_2^lW_{\exp(h_0)})_{2323}(0) =-\kappa\frac{n-3}{2(n-2)}(l+2)!\ne0.
\end{gathered}
\end{equation}
Choose the amplitudes small enough to apply Lemma~\ref{lem:metric-extension}. Support separation and independence of $x_1$ give, for $|y|<\rho_N$,
\begin{equation}\label{eq:high-core}
h_0(x_N+y)=a_NH(y/\lambda_N)+\kappa T_\glo^{[l]}(y).
\end{equation}

\par\medskip\noindent\textbf{The projected family.} We treat the additional homogeneous background as an error  relative to the localized copies of $H$.

The additional background is of higher order because $l+2\ge l_0+3$. In all three dimension ranges, $m>p(l_0+2+1/4)$, and the radial estimate \eqref{eq:radial-master} gives

\[
b_N:=\lambda_N^{\min\{l+2,m\}}(1+|\log\lambda_N|) +(\lambda_N/\rho_N)^m=o(a_N).
\]

Set $B_2=\chi a_NH_{ij}(y/\lambda_N)\partial_{ij}u$, and let $w$ be $\mu_N$ times the endpoint linear correction. Linearity and Proposition~\ref{prop:w-decay} give $\|w\|_{C^1(\overline\Omega,\dot H^1(\R^n))}\le Ca_N$. Repeating the estimates of Proposition~\ref{prop:bubble-estimates} and Lemmas~\ref{lem:linear-correction-error}, \ref{lem:compact-sources} for the two terms in \eqref{eq:high-core} yields
\begin{equation}\label{eq:high-source-estimates}
\begin{split}
\|E_{\widehat g}(\chi u)\|_{C^1(\overline\Omega,L^q(M))}&\le C(a_N+b_N),\\
\|E_{\widehat g}(\chi u)-B_2\|_{C^1(\overline\Omega,L^q(M))}&\le C(a_N^p+b_N),\\
\|\LL_{\widehat g,\chi u}(\chi w)+B_2\|_{C^1(\overline\Omega,L^q(M))} &\le C(a_N^p+b_N).
\end{split}
\end{equation}
Indeed, the quadratic endpoint and mixed sources use \eqref{eq:radial-master} with exponents $2(l_0+2)$ and $l_0+l+4$ and coefficients $\mu_N^2$ and $\mu_N$, respectively. Their inner contributions are $O(a_N^2+a_N\lambda_N^{l+2})$; the exterior radial contributions, including borderline logarithms, are $O(a_N^p)$ by  $m>p(l_0+2+1/4)$. The background sources use exponents $l+2,2(l+2)$ and are $O(b_N)$. The cutoff errors are absorbed in $b_N$. The constraint defects have $C^1(\overline\Omega)$ norm bounded by $C\mu_N\lambda_N^n$, which is also absorbed in $b_N$. Since $\mu_N$ and $\kappa$ are independent of the normalized parameters, the same estimates apply to the first parameter derivatives.

The rescaled metric tends to $\delta$ and the differentiated error term is $o(1)$, so the proof of Proposition~\ref{prop:compact-inverse} applies. Proposition~\ref{prop:compact-projected-family}, with the source bounds \eqref{eq:high-source-estimates}, then gives
\begin{equation}\label{eq:high-corrections}
\begin{aligned}
\|v-\chi u\|_{C^1(\overline\Omega,H^1(M))}&\le Ca_N,\\
\|v-\chi u-\chi w\|_{C^1(\overline\Omega,H^1(M))} &\le C(a_N^p+b_N)=o(a_N).
\end{aligned}
\end{equation}
\par\medskip\noindent\textbf{Reduced-functional comparison.} We compare the reduced functional with $a_N^2F$; we estimate the correction error and the additional geometric terms in $C^1(\overline\Omega)$ and show that their norms, divided by $a_N^2$, tend to zero.

Equation~\eqref{eq:compact-energy-exact} gives
\begin{equation}\label{eq:high-correction-energy}
\begin{aligned}
&\Bigl\|J_{\widehat g}^M(v)-J_{\widehat g}^M(\chi u)-\int_M B_2\chi w\,dV_{\widehat g}\Bigr\|_{C^1(\overline\Omega)}\\
&\qquad\le C\bigl[(a_N^p+b_N)(a_N+b_N)+(a_N+b_N)^{p+1}\bigr]=o(a_N^2).
\end{aligned}
\end{equation}

It remains to check  the terms involving the new background. The first-order terms cancel by \eqref{eq:global-stress}, and the endpoint quadratic term is $a_N^2F$ by the scaling and replacement argument of Lemma~\ref{lem:polynomial-replacement}. As in Lemma~\ref{lem:metric-energy}, put $K_A=\lambda_N^{n-2}\int_0^{\rho_N} (\lambda_N+r)^{A+2-2n}r^{n-1}\,dr$. Splitting at $r=\lambda_N$ gives $K_A\le C(\lambda_N^A+\lambda_N^{n-2})(1+|\log\lambda_N|)$. The pure background and mixed quadratic terms are bounded by $C(K_{2(l+2)}+\mu_NK_{l_0+l+4})$, and Young's inequality bounds the cubic terms by $C(\mu_N^3K_{3(l_0+2)}+K_{3(l+2)})$. All are $o(a_N^2)$: the relevant quadratic powers satisfy
\[
2(l-l_0)- 1/2\ge 3/2,\quad l-l_0-1/4\ge 3/4,\quad n-2-2(l_0+2)-1/2>0.
\]
while the inner endpoint cubic term gains $a_N$ and the inner background cubic term gains the positive power $3(l+2)-2(l_0+2)-1/2$. These strict gaps also absorb any logarithm. The cutoff and exterior errors are $O((\lambda_N/\rho_N)^{n-2})=o(a_N^2)$ by the same tail estimates as in \eqref{eq:cutoff-correction-energy} and Lemma~\ref{lem:polynomial-replacement}. Normalized derivatives have the same radial majorants, so  the same estimates hold in $C^1(\overline\Omega)$. Together with \eqref{eq:high-correction-energy}, they give
\[
\|\widehat F_N-a_N^2F\|_{C^1(\overline\Omega)}=o(a_N^2),
\]
which proves \eqref{eq:higher-order-energy-limit}.
\end{proof}


\begin{thebibliography}{99}
\bibitem{Almaraz} S.~Almaraz, \emph{Blow-up phenomena for scalar-flat metrics on manifolds with boundary}, J. Differential Equations \textbf{251} (2011), no.~7, 1813--1840.
\bibitem{AmbrosettiMalchiodi} A.~Ambrosetti and A.~Malchiodi, \emph{A multiplicity result for the Yamabe problem on $S^n$}, J. Funct. Anal. \textbf{168} (1999), no.~2, 529--561.
\bibitem{Aubin} T.~Aubin, \emph{\'{E}quations diff\'{e}rentielles non lin\'{e}aires et probl\`eme de Yamabe concernant la courbure scalaire}, J. Math. Pures Appl. (9) \textbf{55} (1976), 269--296.
\bibitem{BertiMalchiodi} M.~Berti and A.~Malchiodi, \emph{Non-compactness and multiplicity results for the Yamabe problem on $S^n$}, J. Funct. Anal. \textbf{180} (2001), no.~1, 210--241.
\bibitem{BiHaoHeShiZhu} Y.~Bi, T.~Hao, S.~He, Y.~Shi, and J.~Zhu, \emph{A proof for the Riemannian positive mass theorem up to dimension 19}, preprint (2026), \href{https://arxiv.org/abs/2603.02769}{arXiv:2603.02769}.
\bibitem{Brendle} S.~Brendle, \emph{Blow-up phenomena for the Yamabe equation}, J. Amer. Math. Soc. \textbf{21} (2008), no.~4, 951--979.
\bibitem{BrendleMarques} S.~Brendle and F.~C.~Marques, \emph{Blow-up phenomena for the Yamabe equation II}, J. Differential Geom. \textbf{81} (2009), no.~2, 225--250.
\bibitem{BrendleWang} S.~Brendle and Y.~Wang, \emph{A dimension descent scheme for the positive mass theorem in arbitrary dimension}, preprint (2026), \href{https://arxiv.org/abs/2604.08473}{arXiv:2604.08473}.
\bibitem{BrezisKato} H.~Brezis and T.~Kato, \emph{Remarks on the Schr\"odinger operator with singular complex potentials}, J. Math. Pures Appl. (9) \textbf{58} (1979), 137--151.
\bibitem{CGS} L.~A.~Caffarelli, B.~Gidas, and J.~Spruck, \emph{Asymptotic symmetry and local behavior of semilinear elliptic equations with critical Sobolev growth}, Comm. Pure Appl. Math. \textbf{42} (1989), no.~3, 271--297.
\bibitem{ChodoshMantoulidisSchulze}
O.~Chodosh, C.~Mantoulidis, and F.~Schulze, \emph{Generic regularity for minimizing hypersurfaces in dimensions 9 and 10}, Publ. Math. Inst. Hautes \'Etudes Sci. \textbf{143} (2026), 143--188.

\bibitem{ChodoshMantoulidisSchulzeWang}
O.~Chodosh, C.~Mantoulidis, F.~Schulze, and Z.~Wang, \emph{Generic regularity for minimizing hypersurfaces in dimension 11}, preprint (2025), \href{https://arxiv.org/abs/2506.12852}{arXiv:2506.12852}.
\bibitem{DengLuWei} B.~Deng, H.~Lu, and J.~Wei, \emph{Blow-up phenomena for the $\sigma_2$-Yamabe equation}, preprint (2026), \href{https://arxiv.org/abs/2607.11067}{arXiv:2607.11067}.
\bibitem{Druet1} O.~Druet, \emph{From one bubble to several bubbles: the low-dimensional case}, J. Differential Geom. \textbf{63} (2003), no.~3, 399--473.
\bibitem{Druet2} O.~Druet, \emph{Compactness for Yamabe metrics in low dimensions}, Int. Math. Res. Not. \textbf{2004} (2004), no.~23, 1143--1191.
\bibitem{GeWangWei} Y.~Ge, G.~Wang, and W.~Wei, \emph{A Yamabe problem for the quotient between the $Q$ curvature and the scalar curvature}, preprint (2026), \href{https://arxiv.org/abs/2603.15074}{arXiv:2603.15074}.
\bibitem{GilbargTrudinger} D.~Gilbarg and N.~S.~Trudinger, \emph{Elliptic Partial Differential Equations of Second Order}, reprint of the 1998 edition, Classics in Mathematics, Springer-Verlag, Berlin, 2001.
\bibitem{GolubitskyGuillemin} M.~Golubitsky and V.~Guillemin, \emph{Stable Mappings and Their Singularities}, Graduate Texts in Mathematics, vol.~14, Springer-Verlag, New York, 1973.
\bibitem{GongKimWei} L.~Gong, S.~Kim, and J.~Wei, \emph{Compactness and non-compactness theorems of the fourth- and sixth-order constant $Q$-curvature problems}, preprint (2025), \href{https://arxiv.org/abs/2502.14237}{arXiv:2502.14237}.
\bibitem{GongLi} L.~Gong and Y.~Li, \emph{Conformal metrics of constant scalar curvature with unbounded volumes}, Proc. Lond. Math. Soc. (3) \textbf{131} (2025), no.~1, Paper No.~e70069, 55~pp.
\bibitem{KMS} M.~A.~Khuri, F.~C.~Marques, and R.~M.~Schoen, \emph{A compactness theorem for the Yamabe problem}, J. Differential Geom. \textbf{81} (2009), no.~1, 143--196.
\bibitem{LeeParker} J.~M.~Lee and T.~H.~Parker, \emph{The Yamabe problem}, Bull. Amer. Math. Soc. (N.S.) \textbf{17} (1987), no.~1, 37--91.
\bibitem{LesourdUngerYau} M.~Lesourd, R.~Unger, and S.-T.~Yau, \emph{The positive mass theorem with arbitrary ends}, J. Differential Geom. \textbf{128} (2024), 257--293.
\bibitem{LiWangWei} C.~Li, G.~Wang, and W.~Wei, \emph{Blow-up phenomena for the constant $Q/R$-curvature equation}, preprint (2026), \href{https://arxiv.org/abs/2604.20571}{arXiv:2604.20571}.
\bibitem{LiXiong} Y.Y.~Li and J.~Xiong, \emph{Compactness of conformal metrics with constant $Q$-curvature. I}, Adv. Math. \textbf{345} (2019), 116--160.
\bibitem{LiZhang1} Y.Y.~Li and L.~Zhang, \emph{A Harnack type inequality for the Yamabe equation in low dimensions}, Calc. Var. Partial Differential Equations \textbf{20} (2004), no.~2, 133--151.
\bibitem{LiZhang2} Y.Y.~Li and L.~Zhang, \emph{Compactness of solutions to the Yamabe problem. II}, Calc. Var. Partial Differential Equations \textbf{24} (2005), 185--237.
\bibitem{LiZhang3} Y.Y.~Li and L.~Zhang, \emph{Compactness of solutions to the Yamabe problem. III}, J. Funct. Anal. \textbf{245} (2007), no.~2, 438--474.
\bibitem{LiZhu} Y.Y.~Li and M.~Zhu, \emph{Yamabe type equations on three-dimensional Riemannian manifolds}, Commun. Contemp. Math. \textbf{1} (1999), no.~1, 1--50.
\bibitem{Lohkamp}
J.~Lohkamp, \emph{The higher dimensional positive mass theorem I}, preprint (2006), revised 2016, \href{https://arxiv.org/abs/math/0608795}{arXiv:math/0608795}.
\bibitem{MarquesApriori} F.~C.~Marques, \emph{A priori estimates for the Yamabe problem in the non-locally conformally flat case}, J. Differential Geom. \textbf{71} (2005), no.~2, 315--346.
\bibitem{MarquesBlowup} F.~C.~Marques, \emph{Blow-up examples for the Yamabe problem}, Calc. Var. Partial Differential Equations \textbf{36} (2009), no.~3, 377--397. \href{https://webdoc.sub.gwdg.de/ebook/serien/e/IMPA_A/546.pdf}{Preprint}.
\bibitem{Nirenberg}
L.~Nirenberg,, \emph{Topics in nonlinear functional analysis}, Courant Lect. Notes Math., 6 New York University, Courant Institute of Mathematical Sciences, New York; American Mathematical Society, Providence, RI, 2001, xii+145 pp.
\bibitem{Schoen} R.~M.~Schoen, \emph{Conformal deformation of a Riemannian metric to constant scalar curvature}, J. Differential Geom. \textbf{20} (1984), no.~2, 479--495.
\bibitem{SchoenLCF} R.~M.~Schoen, \emph{On the number of constant scalar curvature metrics in a conformal class}, in \emph{Differential Geometry}, Pitman Monographs and Surveys in Pure and Applied Mathematics, vol.~52, Longman Scientific \& Technical, 1991, pp.~311--320.
\bibitem{SchoenYauMass} R.~M.~Schoen and S.-T.~Yau, \emph{On the proof of the positive mass conjecture in general relativity}, Comm. Math. Phys. \textbf{65} (1979), 45--76.
\bibitem{SchoenYauSingularities}
R.~M.~Schoen and S.-T.~Yau, \emph{Positive scalar curvature and minimal hypersurface singularities}, Surv. Differ. Geom. \textbf{24} (2019), no.~1, 441--480.
\bibitem{Smale}
N.~Smale, \emph{Generic regularity of homologically area minimizing hypersurfaces in eight-dimensional manifolds}, Comm. Anal. Geom. \textbf{1} (1993), no.~2, 217--228.
\bibitem{Trudinger} N.~S.~Trudinger, \emph{Remarks concerning the conformal deformation of Riemannian structures on compact manifolds}, Ann. Scuola Norm. Sup. Pisa Cl. Sci. (3) \textbf{22} (1968), 265--274.
\bibitem{WeiZhao} J.~Wei and C.~Zhao, \emph{Non-compactness of the prescribed $Q$-curvature problem in large dimensions}, Calc. Var. Partial Differential Equations \textbf{46} (2013), no.~1--2, 123--164.
\bibitem{Witten} E.~Witten, \emph{A new proof of the positive energy theorem}, Comm. Math. Phys. \textbf{80} (1981), 381--402.
\bibitem{Yamabe} H.~Yamabe, \emph{On a deformation of Riemannian structures on compact manifolds}, Osaka Math. J. \textbf{12} (1960), no.~1, 21--37.
\end{thebibliography}
\end{document}